\documentclass[11pt]{article}
\usepackage{caption}
\usepackage{appendix}
\usepackage{amsthm}
\usepackage{multirow}
\usepackage{hyperref}
\usepackage{setspace}
\usepackage[utf8]{inputenc}
\usepackage{graphicx}
\usepackage{tikz}
\usetikzlibrary{bayesnet, backgrounds, positioning}
\tikzset{
	treenode/.style = {shape=rectangle, rounded corners,
		draw, align=center,
		top color=white, bottom color=blue!20},
	root/.style     = {treenode, font=\Large, bottom color=red!30},
	env/.style      = {treenode, font=\ttfamily\normalsize},
	dummy/.style    = {circle,draw}
}	
	\newcommand{\blind}{0}
	
    \makeatletter
    \renewcommand\section{\@startsection {section}{1}{\z@}%
                                       {-3.5ex \@plus -1ex \@minus -.2ex}%
                                       {2.3ex \@plus.2ex}%
                                       {\normalfont\fontfamily{phv}\fontsize{16}{19}\bfseries}}
    \renewcommand\subsection{\@startsection{subsection}{2}{\z@}%
                                         {-3.25ex\@plus -1ex \@minus -.2ex}%
                                         {1.5ex \@plus .2ex}%
                                         {\normalfont\fontfamily{phv}\fontsize{14}{17}\bfseries}}
    \renewcommand\subsubsection{\@startsection{subsubsection}{3}{\z@}%
                                        {-3.25ex\@plus -1ex \@minus -.2ex}%
                                         {1.5ex \@plus .2ex}%
                                         {\normalfont\normalsize\fontfamily{phv}\fontsize{14}{17}\selectfont}}
    \makeatother
	
	\usepackage{amsmath}
	\usepackage{graphicx}
	\usepackage{enumerate}
	\usepackage{natbib} 
	\usepackage{url} 

\usepackage{xargs}                      
\usepackage[pdftex,dvipsnames]{xcolor}  
\usepackage[colorinlistoftodos,prependcaption,textsize=tiny]{todonotes}
\newcommandx{\unsure}[2][1=]{\todo[linecolor=red,backgroundcolor=red!25,bordercolor=red,#1]{#2}}
\newcommandx{\change}[2][1=]{\todo[linecolor=blue,backgroundcolor=blue!25,bordercolor=blue,#1]{#2}}
\newcommandx{\info}[2][1=]{\todo[linecolor=OliveGreen,backgroundcolor=OliveGreen!25,bordercolor=OliveGreen,#1]{#2}}
\newcommandx{\improvement}[2][1=]{\todo[linecolor=Plum,backgroundcolor=Plum!25,bordercolor=Plum,#1]{#2}}
\newcommandx{\thiswillnotshow}[2][1=]{\todo[disable,#1]{#2}}
\usepackage{comment}
 \usepackage{todonotes} 

\usepackage{subfig}

 \newtheorem{lemma}{Lemma}

\newtheorem{prop}{Proposition}

\newtheorem{corollary}{Corollary}

\usepackage{MnSymbol}%
\usepackage{wasysym}%
\usepackage{mathtools}
\usepackage{amsmath}
\usepackage{algorithm}
\usepackage{algpseudocode}
\usepackage{tablefootnote}

\begin{document}
		
		\def\spacingset#1{\renewcommand{\baselinestretch}%
			{#1}\small\normalsize} \spacingset{1}
		
		\if0\blind
		{
			\title{\bf Fortifying Distribution Network Nodes Subject to Network-Based Disruptions}
			\author{Pelin Ke\c{s}rit\thanks{pelin@tamu.edu} , Bahar \c{C}avdar\thanks{cavdab2@rpi.edu} , and Joseph Geunes\thanks{geunes@tamu.edu}  \\
			$^*$ $^\ddagger$\small{Department of Industrial and Systems Engineering, Texas A\&M University, College Station, TX }\\
             $^\dagger$\small{Department of Industrial and Systems Engineering, Rensselaer Polytechnic Institute, Troy, NY}
             }
			\date{}
			\maketitle
		} \fi
		
		\if1\blind
		{

            \title{\bf \emph{IISE Transactions} \LaTeX \ Template}
			\author{Author information is purposely removed for double-blind review}
			
\bigskip
			\bigskip
			\bigskip
			\begin{center}
				{\LARGE\bf \emph{IISE Transactions} \LaTeX \ Template}
			\end{center}
			\medskip
		} \fi
		\bigskip
		
	\begin{abstract}
 We consider a distribution network for delivering a natural resource or physical good to a set of nodes, each of which serves a set of customers, in which disruptions may occur at one or more nodes. Each node receives flow through a path from a source node, implying that the service at a node is interrupted if one or more nodes on the path from a source node experience a disruption. 
 All network nodes are vulnerable to a future disturbance due to a potential natural or man-made disaster, the severity of which follows some measurable probability distribution.  For each node in the network, we wish to determine a fortification level that enables the node to withstand a disturbance up to a given severity level, while minimizing the expected number of customers who experience a service interruption under a limited fortification budget. We formulate this problem as a continuous, nonlinear knapsack problem with precedence constraints, demonstrate that this optimization problem is $\mathcal{NP}$-Hard for general tree networks and general disturbance severity distributions, and provide a polynomial-time solution algorithm for serial systems, which forms the basis for an effective heuristic approach to problems on tree networks. Our computational test results demonstrate the ability of the proposed heuristic methods to quickly find near-optimal solutions.

	\end{abstract}
			
	\noindent%

	\spacingset{1.5} 

\setlength{\parindent}{15pt}

\section{Introduction} \label{s:intro}

\emph{System Resilience} is often characterized by an ability to recover from disruptions as a result of effective preparation \citep{datta2017supply} and can be thought of as the opposite of vulnerability and brittleness \citep{sharkey2021search}. Improving a system's resilience entails efforts to reduce vulnerability and mitigate uncertainty. As one of the four concepts associated with resilience, enhancing the \emph{robustness} of a system allows for an effective response to random perturbations \citep{woods2015four}. In the context of distribution networks, including critical infrastructure for water, power, gas, and communication, which are essential for everyday life \citep{karakoc2019community}, random perturbations arise in the form of unexpected events such as natural disasters that can reduce or eliminate the ability to function at one or more nodes.  Effectively accounting for the possibility of such disruptions directly impacts a system's reliability and performance, and is therefore of crucial importance in network planning.  This is especially true for power systems \citep{jose2021optimal}, which enable other critical infrastructure systems to operate, including communications, transportation, and health systems \citep{qi2024network}. The performance of a network after a disruption often depends on prior investments in system robustness in addition to  the severity of the perturbation itself \citep{sharkey2021search}. These investments, as part of the preparation phase of resilience planning efforts, include strengthening the underlying network to withstand the negative impacts of a disturbance \citep{biswas2024review}. Such prevention and mitigation strategies in the preparation phase are vital, as they can reduce the resources required to restore the network in the later stages of disaster management \citep{gupta2016disaster}.  Minimizing the negative service impacts that may result from a disturbance requires targeted investment in critical infrastructure components, as a limited fortification budget renders absolute protection infeasible \citep{freiria2015multiscale}. 

 We consider a distribution network for delivering a natural resource, service, or physical good to a set of nodes, each of which serves a set of customers, in which flow disruptions may occur at one or more nodes. Each node receives flow via a path from a source node. If each node receives flow via a single path from a single source, the underlying network contains a tree structure, while the ability for a node to receive flow on multiple paths from the same source introduces network cycles. In either case, a node's ability to provide service to customers depends on the availability of flow from upstream nodes.

Our modeling approach assumes that all network nodes are subject to a common disturbance severity distribution, as may be the case when nodes lie in reasonably close proximity to one another.  We also assume the degree of severity is measurable and follows some well-defined probability distribution. For each node in the network with no cycles, we wish to determine a fortification level that enables the node to withstand a disturbance up to some level of severity, given a limited fortification budget. The problem then seeks to minimize the expected number of customers without service, which is equivalent to maximizing the expected number of customers with no service disruption. We formulate this problem as a continuous, nonlinear knapsack problem with precedence constraints and, although we show that the problem on a general tree is $\mathcal{NP}$-Hard, we characterize the properties of optimal solutions that enable solving the problem for systems with serial networks in polynomial time under mild conditions on the disruption severity distribution.  
  
The difficulty of the resulting optimization problem depends both on the structure of the network and the corresponding disturbance severity distribution. If, for example, the cumulative distribution function (CDF) of the disturbance distribution is concave, then the resulting problem is a convex program and, therefore, can be solved efficiently. Unfortunately, a limited number of continuous probability distributions possess a concave CDF and confining the analysis to this class of distributions severely limits the applicability of our model in practice. For most commonly applied continuous distributions of natural phenomena, such as normal, beta, gamma, logistic, and triangular, the CDF is neither convex nor concave in general. In particular, the CDF often takes the form of a non-decreasing S-shaped curve, which we later describe in greater detail.  For such cases, we develop an algorithm based on the problem's \emph{necessary} (but not sufficient) KKT conditions, which is able to eliminate the computational burden of a global optimization solver that often is unable to find an optimal solution in an acceptable amount of time as the problem size increases. The resulting algorithm is able to solve problems with a serial network structure in polynomial time under mild disturbance distribution assumptions. We use the properties of optimal solutions for series systems to develop a heuristic solution approach for general tree systems. We also demonstrate the effectiveness of using the smallest concave envelope of the CDF to provide strong upper bounds and corresponding heuristic solutions.  

The rest of this paper is organized as follows. Section~\ref{s:litrew} first discusses the related literature, after which Section~\ref{s:probdef} describes the problem setting, introduces our mathematical model, and characterizes useful properties of optimal solutions.  Section~\ref{sec:poly} provides polynomial-time solution algorithms for two special cases: (i) series systems with S-shaped distribution functions, and (ii) general tree systems under a uniform disturbance severity distribution. Motivated by effective methods for solving these special cases, we provide a heuristic solution algorithm for general tree structures and S-shaped distribution functions in Section~\ref{sec:heuristic}.  Section~\ref{s:computational} presents the results of a set of computational study that demonstrate the advantages of our heuristic solution methods over a commercial global optimization solver.  The paper provides concluding remarks and a discussion of future research directions in Section~\ref{s:conclusion}.



\section{Literature Review} \label{s:litrew}

Our proposed model considers fortifying nodes in a distribution network in order to minimize the expected number of affected individuals in the event of a disruption of an uncertain severity under a limited budget. A related stream of (infrastructure) network resilience literature focuses on the preparedness stage before any disruption in an effort to mitigate the impacts of a disturbance through fortification, strengthening, and/or hardening, using a variety of different solution approaches. As we later discuss, our model formulation corresponds to a nonlinear, stochastic, and continuous version of the Knapsack problem with precedence constraints. Therefore, this section reviews both the network resilience literature (Section~\ref{s:litrew:network}) and related literature on knapsack problems (Section~\ref{s:litrew:knapsack}).

\subsection{Network Resilience} \label{s:litrew:network}

Network resilience has been gaining attention in the literature, especially in the context of infrastructure \citep{liu2022resilience}, road \citep{rivera2022road} and power \citep{bhusal2020power} networks. Researchers have employed a variety of methodologies to improve network resilience through fortification under limited resources. In their literature survey, \cite{hoyos2015or} note that disaster operations management studies most commonly use mathematical models with stochastic features. \cite{peeta2010pre} formulate a pre-disaster highway network investment problem as a two-stage stochastic program whose approximation yields a knapsack problem. In a similar problem setting, \cite{yucel2018improving} model link failures using Bayesian dependency and use a tabu search algorithm to determine which links to fortify. In order to maximize freight transportation network resilience with a given budget, \cite{miller2012measuring} develop a two-stage stochastic program and observe that the highest levels or resilience are attained when both preparedness and recovery options are available. \cite{patron2017optimal} strengthen a connected network by increasing the degree of certain nodes subject to attack in a general random network based on each node's degree-to-cost ratio. \cite{banerjee2018hardening} show that determining the nodes to be fortified to minimize the damage in a critical infrastructure network with dependencies is $\mathcal{NP}$-Hard and propose a heuristic solution. \cite{almaleh2022novel} study a microgrid placement problem to ensure the continuity of power at critical nodes in the event of a disaster while minimizing the cost, and show that the corresponding problem is $\mathcal{NP}$-Hard. \cite{costa2024modelling} propose an iterative cutting plane algorithm for fortifying or immunizing an electric transmission network at minimum cost. To find an optimal fortification strategy in decentralized supply systems, \cite{zhang2017fortification} develop a tree-search-based algorithm.

 Although a number of studies focus on improving network resilience or disaster management planning through strengthening efforts under limited resources, the contexts considered and underlying model assumptions vary. \cite{liu2022network} summarize related past work in the context of social, biological, ecological, and infrastructure networks. Most of these studies limit such decisions to whether or not a node is fortified at a predetermined, fixed level, without quantifying the degree of fortification \citep{dehghani2016resource}. Similarly, prior literature often assumes a fixed probability of a disruption \citep{gao2017resilience}. Our work develops an algorithm to fortify nodes in a general tree network against an uncertain disturbance severity level, where the properties of the underlying probability distribution and the network structure determine the problem's complexity.  


\subsection{Knapsack with Precedence Constraints} \label{s:litrew:knapsack}

The precedence-constrained knapsack problem (PCKP) is a generalization of the 0-1 knapsack problem (KP), where an item may be included in the knapsack only if all of its predecessors are included \citep{wilbaut2022knapsack}. A variety of practical problems can be modeled as a PCKP, including capital budgeting \citep{salkin1975knapsack}, open-pit mine scheduling \citep{moreno2010large, samavati2017methodology, nancel2022recursive}, machine scheduling \citep{wikum1994one}, production scheduling \citep{bienstock2010solving}, and subject selection in personalized learning \citep{aslan2023precedence}.

The PCKP is $\mathcal{NP}$-Hard \citep{garey1979computers} which can be shown by a reduction from the KP. \cite{ibarra1978approximation} develop polynomial-time approximation algorithms for special cases of the PCKP where the degree of the polynomial depends on the accuracy level of the approximate solution.
\cite{johnson1983knapsacks} use dynamic programming to develop pseudopolynomial time 
algorithms and fully polynomial-time approximation schemes for the PCKP. \cite{boyd1993polyhedral} shows that the convex hull of feasible points for the precedence constraints (without the knapsack constraint) is the same as its LP relaxation feasible region and characterizes facets for variations of the problem. \cite{wilbaut2022knapsack} note that this enables solving the Lagrangian relaxation of the PCKP (where the knapsack constraint is relaxed) in polynomial time, although the integrality property implies that the resulting Lagrangian relaxation bound will be the same as that of the PCKP LP relaxation value \citep[see, e.g.,][]{geoffrion2009lagrangean}. \cite{you2007pegging} use a pegging test for the LP relaxation of the PCKP to reduce the problem size by removing variables that must take a value of either 0 or 1. 

Related literature includes a plethora of studies that tackle the PCKP using a variety of methods that either present exact solutions to the problem under simplifying assumptions or provide heuristic solutions. \cite{park1997lifting} use lifting cover inequalities for cutting plane algorithms for the solution of problems that have the PCKP as a substructure. Similarly, \cite{van1999lifting} study valid inequalities derived from minimal induced covers to identify a class of lifting coefficients for the PCKP and show that lifting coefficients can be obtained in polynomial time for the tree knapsack problem, which is a special case of the PCKP where the precedence graph is a tree. \cite{boland2012clique} show that adding clique-based inequalities to the PCKP helps to decrease the LP relaxation gap. Representing the PCKP as a directed acyclic graph, \cite{samphaiboon2000heuristic} adopt a dynamic programming approach to find an optimal solution to the PCKP for small instances, and develop heuristic algorithms based on preprocessing techniques for larger instances. \cite{samavati2017methodology} develop an LP-based heuristic algorithm by adding cover and clique constraints for the PCKP. \cite{nancel2022recursive} use a recursive time aggregation-disaggregation heuristic after reducing the size of the PCKP.


The problem of determining the fortification levels in a distribution network subject to disruptions can be viewed as a version of the knapsack problem as the problem considers allocating a limited budget to nodes whose weights are implicitly determined by the number of customers they serve. Each node receives flow through a single path from the source; therefore, for a node to be fortified at a certain level, every node along that path from the source also needs to be fortified. This property introduces precedence constraints to the problem, leading to a PCKP. However, because the value of fortifying a node depends on an uncertain disturbance severity level and the probability distribution it follows, this results in a nonlinear, stochastic, and continuous generalization of the PCKP, which requires employing a solution approach that differs from the methodologies commonly applied in the aforementioned literature.


\section{Problem Definition and Model Properties} \label{s:probdef}
We first define the problem and provide a general formulation of our model in Section \ref{sec:PF}. Section \ref{sec:MPR} then establishes key properties of optimal solutions that permit reformulating the problem as a nonlinear stochastic knapsack problem with precedence constraints.  


\subsection{Problem Formulation}\label{sec:PF} 
Consider a commodity flow tree network $G(N,A)$, where each node $i \in N$ serves a given integer number of customers, $w_i$, using the flow of a single commodity that originates at a source node denoted as node $1$. 
Flow is distributed to each node via a path consisting of its predecessors. Therefore, in order for node $i$ to be able to serve its customers, every node on the path connecting node $i$ to the source node (including node $i$ itself) denoted by $\mathcal{P}(i)$, must function properly (i.e., each node in $\mathcal{P}(i)$ must not have failed). 
We assume that all nodes are subject to a set of adverse and uncertain future events and that we can normalize the maximum severity of these events, which we denote as the random variable $Y$, on a scale between 0 and 1, with 1 corresponding to the most severe and 0 the least severe. Examples of severity measures include maximum wind speed, highest flood level, or earthquake intensity. 
We also assume that each node $i$ can be fortified to withstand a severity level up to $x_i$ at a cost $c_ix_i$ for each $i\in N$, where $c_i$ is the per-unit cost of fortification. The fortification level is measured using the same scale as the disruption level and, therefore, takes a value in $[0,1]$. 
Thus, a node fails if the severity of the worst event in a season ($Y$) is greater than or equal to the node's fortification level. That is, node $i$ fails if $Y \geq x_i$.
Then, the probability 
that the service at node $i$ is disrupted (i.e., that node $i$ is without servive) is written as
\[P\left(Y \geq \min_{j\in \mathcal{P}(i)}\{x_j\}\right).\]
This is the probability that the severity of disruption exceeds the fortification level of the weakest node on path $\mathcal{P}(i)$. We assume that all nodes initially have a base level of fortification, which can be normalized to zero.
\noindent Table~\ref{t.param} presents a summary of our notation.



\begin{table}[h]
    \caption{Notation} 
    \centering
    \begin{tabular}{|rl|}
        \hline
        \multicolumn{2}{|l|}{\textbf{Parameters}}\\
        \hline
        $c_i $ & Per-unit fortification cost at node $i$\\
        $\mathcal{P}(i)$ & Set of nodes on the path from source to node $i$\\
        $\pi(i)$ & The immediate predecessor of node $i$\\
        $\sigma(i)$ & Set of immediate successors of node $i$ \\
        $\mathcal{L}$ & Set of leaf nodes \\
        $w_i$ & Number of customers solely served at node $i$\\
        $B$ & Total fortification budget\\
        $n$ & Number of nodes, equal to $|N|$\\
        \hline
        \multicolumn{2}{|l|}{\textbf{Decision Variables}}\\
        \hline
        $x_i$ & Fortification level of node $i$\\
        \hline
    \end{tabular}
    \label{t.param}
\end{table} 

We wish to determine the node fortification levels ($x_i$ values) under a limited budget $B$ that minimize the expected number of customers without service after a disruption.  We formulate this problem 
as follows:
\begin{align}
    \text{Minimize} \quad \ & \sum_{i=1}^n P\left(Y \geq \min_{j\in \mathcal{P}(i)}\{x_j\}\right)w_i \label{m.1}\\
    \text{Subject to: } & \sum_{i=1}^n c_ix_i \leq B \label{m.2}\\
    & 0 \leq x_i \leq 1, \hspace{3mm}  i = 1,\ldots,n. \label{m.3} 
\end{align}

The objective function (\ref{m.1}) minimizes the expected number of customers without service, which is the sum of the probability that each node loses service multiplied by the number of customers solely served by that node. Constraint (\ref{m.2}) ensures that the total cost of fortification does not exceed the fortification budget, $B$. Constraints (\ref{m.3}) ensure that the fortification level of each node takes a value on the $[0,1]$ interval.
To avoid trivial cases where all nodes can be fully fortified, i.e., $x_i=1$ for all $i\in N$, we assume that $B < \sum_{i=1}^{n} c_i$.

\subsection{Model Properties and Reformulation}\label{sec:MPR}

The minimization objective function (\ref{m.1}) implies the following proposition (all proofs can be found in Appendix \ref{sec:Proofs}).

\begin{prop} \label{p1}
    In an optimal solution for the network fortification problem, $x_i \geq x_j$ for all $i \in \mathcal{P}(j)$ where node $i$ is a predecessor of node $j$, for any $c_j$, $j\in \{1,\dots,n\}$.
\end{prop}

\begin{corollary}  \label{c1}
In an optimal solution for the network fortification problem, if $x_i=0$, then $x_j=0$ for each successor node $j$ of node $i$ in the underlying infrastructure network.
\end{corollary}

As noted in Proposition~\ref{p1}, $x_j \geq x_i$ for all $j \in \mathcal{P}(i)$ where node $j$ is the predecessor of node $i$. Hence, it follows that 
\begin{equation}
    P(y \geq \min_{j\in \mathcal{P}(i)}\{x_j\})=P(y \geq x_i)=1-P(y \leq x_i) = 1-F_Y(x_i), \quad \forall i \in N.
\end{equation}
The objective function (\ref{m.1}) can thus be rewritten as $\sum_{i=1}^n [1-F_Y(x_i)]w_i$. That is, minimizing the expected number of customers without service is equivalent to maximizing the expected number of customers without a service disruption.



\begin{figure}[h]
	\centering
    \resizebox{60mm}{60mm}{
	\begin{tikzpicture}[thick]      
    \node [latent] (source) {1};
	\node [latent, yshift=-3mm, below=of source] (n2) {2};
    \node [const, right=of source, xshift=-7mm] (stext) {Source};
	\node [latent, left=of n2, yshift=3mm] (n1) {3};
	\node [latent, left=of n1] (n4) {4};
	\node [latent, right=of n2, yshift=3mm] (n3) {5};
	\node [latent, right=of n3, yshift=7mm] (n5) {6};
	\node [latent, right=of n3, yshift=-7mm] (n6) {7};
	\node [latent, below=of n2] (n7) {8};
	\node [latent, below=of n7] (n8) {9};
	\node [latent, left=of n8] (n9) {10};
	\node [latent, below=of n8] (n10) {11};
	\node [latent, right=of n8] (n11) {12};
	\node [latent, right=of n11] (n12) {13};
	\node [latent, below=of n12] (n13) {14};
	\edge {source}{n1,n2,n3};
	\edge {n1}{n4};
	\edge {n3}{n5,n6};
	\edge {n2}{n7};
        \edge {n7}{n8};
        \edge {n8}{n9,n10,n11};
        \edge {n11}{n12};
        \edge {n12}{n13};
	\end{tikzpicture}
    }
	\caption{A general tree topology for a distribution network.}\label{fig:treegraph}
\end{figure}
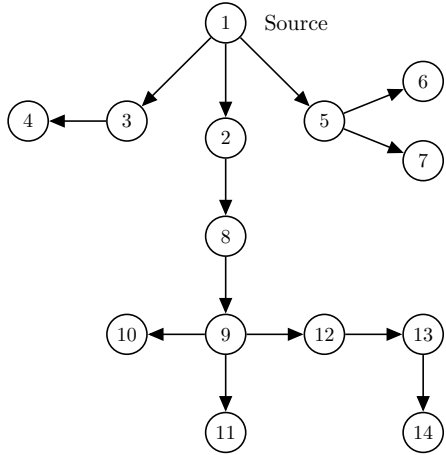

Using Proposition~\ref{p1}, the mathematical model characterized by (\ref{m.1})-(\ref{m.3}) can be rewritten in the following form for a distribution network with no cycles (which will be referred to as a general tree, please see Figure~\ref{fig:treegraph} for an example):
\begin{eqnarray}
\mbox{[$KP_F$]} \quad    \text{Maximize} \ \ & \sum_{i=1}^n F_Y(x_i)w_i & \hspace{4cm} \label{gt_obj}\\
    \text{Subject to: } & \sum_{i=1}^n c_ix_i \leq B,  & \label{eq:capacity}\\
    & x_1 \leq 1, & \label{eq:x1}\\
    & x_i - x_{\pi(i)} \leq 0,   & i \in N \backslash \{1\}, \label{eq:xi} \\
    & x_i \geq 0, & i \in N. \label{gt_mod}
\end{eqnarray}

In the special case that $F_Y(y)$ is a concave function, then because a weighted sum of concave functions is also concave, [$KP_F$] is a convex program. Thus, the Karush Kuhn-Tucker (KKT) optimality conditions, provided in Section \ref{sec:KKT}, are necessary and sufficient for optimality \citep[see, e.g.,][]{bazaraa2006nonlinear}. Distributions with concave CDFs (e.g., exponential, Pareto) therefore yield convex programs and can be solved efficiently using a commercial solver. If $F_Y(y)$ is concave for all $y$ greater than or equal to the distribution's mode (which in not uncommon, as discussed below), then if the minimum fortification level equals the mode, we obtain an objective function that consists of the sum of concave functions, which also yields a convex program. Conversely, if $F_Y(y)$ is convex for $y$ less than or equal to the mode and the maximum fortification level equals the mode, we obtain a convex objective function, implying an extreme point optimal solution.

For a number of broadly applicable distributions (e.g., Cauchy, gamma, Erlang, log-normal, logistic, normal, triangular, and Weibull)\footnote{Many commonly applied distributions (e.g., normal, logistic, gamma) have unbounded supports, and thus, for use in our model, their parameters need to be chosen such that all but a negligible part of the distribution values fall on the $[0,1]$ interval $[a,b]$. For example, for a normal distribution, choosing a mean $\mu$ and standard deviation $\sigma$ such that $\sigma \le \min\left\{\frac{\mu}{3},\frac{1-\mu}{3}  \right\}$ ensures that more than $99.7\%$ of the distribution's values fall on the $[0,1]$ interval, and a negligible part of the distribution falls outside these bounds. For distributions with bounded supports (uniform, triangular, beta), on the other hand, it is possible to choose upper and lower bounds within the desired range.},  $F_Y(y)$ is neither convex nor concave in general, but is convex and nondecreasing for values less than the mode, which we denote as $\beta$, and concave and nondecreasing thereafter, i.e., $F_Y(y)$ is convex for $y\leq \beta$ and is concave for $y\geq \beta$, where $a=0 < \beta < b=1$. Such functions are categorized as sigmoid functions and are often referred to as S-shaped or simply S-curves.  \cite{augrali2009solving} show that the continuous knapsack problem with S-curve return functions but without precedence constraints is $\mathcal{NP}$-Hard when the return functions associated with different items are not identical. Despite the fact that the distribution function in problem $[KP_F]$ is identical for each node, we are able to show that the following proposition holds.

\begin{prop} \label{p:nphard}
Problem $[KP_F]$ is $\mathcal{NP}$-Hard.
\end{prop}

Proposition \ref{p:nphard} implies that for general tree structures and general probability distributions, the time required to solve $[KP_F]$ grows exponentially in the problem size (unless $\mathcal{P}=\mathcal{NP}$).  Thus, heuristic solution methods will be required to find effective solutions in an acceptable computing time. Section \ref{sec:series} proposes a set of heuristic solution methods motivated primarily by the ability to solve the special case of $[KP_F]$ on a serial network in polynomial time under mild disturbance severity distribution assumptions. The ability to do this relies on properties resulting from analysis of the problem's KKT conditions, provided in the next subsection.
\vspace{-6mm}
\subsection{KKT Conditions}\label{sec:KKT}
When $F_Y(y)$ is continuous and differentiable but non-concave, then because all of the constraints are linear, the KKT conditions are \emph{necessary}, but not sufficient, for optimality.  We next analyze these conditions and use them to establish an important property of KKT points that will be very useful in both deriving a polynomial-time solution for series systems (i.e., the degree of each node is at most two and the network is connected), as well as in the development of heuristic solution methods for general tree systems.

To write the KKT conditions for $[KP_F]$, denote $\alpha$ as the KKT multiplier associated with Constraint \ref{eq:capacity} and $\lambda_1$ the KKT multiplier associated with Constraint \ref{eq:x1}.  Similarly, $\lambda_i$ denotes the KKT multiplier associated with Constraint \ref{eq:xi}, for $i \in N\backslash\{1\}$. The KKT conditions are as follows.

\noindent \textbf{KKT Conditions for} $\mathbf{[KP_F]:}$ 
\begin{eqnarray}
    & [KP_F] \ \eqref{eq:capacity} - \eqref{gt_mod}, & \nonumber \\
    &w_i\frac{\partial F_Y(x_i)}{\partial x_i}=\alpha c_i + \lambda_i -\sum_{j\in \sigma(i)}\lambda_{j}, \hspace{1cm} & i \in N \backslash \mathcal{L}, \label{eq:kkt5} \\
    &w_{i} \frac{\partial F_Y(x_{i})}{\partial x_{i}} \leq \alpha c_{i} + \lambda_{i}, & i \in \mathcal{L}, \label{eq:kkt6} \\
    &\alpha(B-\sum_{i=1}^n c_ix_i)=0,& \label{eq:kkt7}\\
    &\lambda_1(1-x_1)=0,&  \label{eq:kkt8} \\
    &\lambda_i(x_{\pi(i)}-x_i)=0, & i \in N \backslash \{1\}  \label{eq:kkt9} \\
    &\alpha \geq 0,& \label{eq:kkt10}\\
    &\lambda_i \geq 0, & i\in N. \label{eq:kkt11}
\end{eqnarray}
In addition to the initial $[KP_F]$ primal feasibility conditions, conditions \eqref{eq:kkt5}-\eqref{eq:kkt6} correspond to dual feasibility constraints, while \eqref{eq:kkt7}-\eqref{eq:kkt9} represent complementary slackness conditions.  The KKT multiplier nonnegativity conditions \eqref{eq:kkt10}-\eqref{eq:kkt11} complete the set of conditions required for a KKT point.


Next, suppose we have a set of nodes $I$ such that $x_i=x_j$ for all $(i,j)\in I$. Let $I_s \subseteq I$ denote the set of $i \in I$ such that $\pi(i) \notin I$, that is those whose predecessor is not in $I$. Let $I_p \subseteq I$ denote the set of $i \in I$ such that $\sigma(i) \cap I= \emptyset$ and $\pi(i) \in I$, i.e., those whose predecessor is in $I$ but which have no successors in $I$. Because $x_i<x_{\pi(i)}$ for $i \in I_s$, condition \eqref{eq:kkt9} implies that $\lambda_i=0$ for $i \in I_s$. In addition, for $i \in I_p$, $x_i > x_j$ for all $j\in \sigma(i)$, which implies $\lambda_j=0$ for all $j \in \pi(i)$ when $i \in I_p$. 
Therefore, for the nodes $i\in I$, we can rewrite the condition in Equation \ref{eq:kkt5} as follows
\begin{align*}
    & w_i\frac{\partial F_Y(x_i)}{\partial x_i}=\alpha c_i - \sum_{j\in \sigma(i) \cap I}\lambda_{j},& i \in I_s, \\
    & w_i\frac{\partial F_Y(x_i)}{\partial x_i}=\alpha c_i +\lambda_i - \sum_{j\in \sigma(i) \cap I}\lambda_{j}, & i \in I \backslash \{I_s \cup I_p\}, \\
    & w_i\frac{\partial F_Y(x_i)}{\partial x_i}=\alpha c_i + \lambda_i, & i \in  I_p.
\end{align*}
Note that our assumption that $x_i=x_j$ for all $(i,j) \in I$ implies that $\frac{\partial F_Y(x_i)}{\partial x_i}=\frac{\partial F_Y(x_j)}{\partial x_j}$ for all $(i,j) \in I$. 
Letting $W(I)=\sum_{i \in I} w_i$ and $C(I)=\sum_{i \in I}c_i$, and summing the above equations over all $i \in I$ gives
\begin{center}
    $\frac{\partial F_Y(x_i)}{\partial x_i}W(I)=\alpha C(I)$
\end{center}
i.e.,
\begin{eqnarray}
    \alpha=\frac{W(I)}{C(I)}f_Y(x_i), \label{eqn:margcont}
\end{eqnarray}
where $f_Y(x)$ corresponds to the probability density function (pdf) of $Y$ at $Y=y$. This implies that at optimality, all subsets of items with equal fortification levels must have identical values of the weighted marginal contribution per unit of capacity consumption, $\left( \frac{W(I)}{C(I)}f_Y(x_i) \right)$, and this must equal the value of the KKT multiplier associated with the budget constraint \eqref{eq:capacity}. 

\section{Polynomially Solvable Special Cases}\label{sec:poly}
This section considers two special cases that permit a polynomial-time solution.  The first case deals with serial systems where the distribution function $F_Y(y)$ takes the form of a differentiable S-curve.  The second case, discussed in Section \ref{sec:uniform}, considers general tree systems under a uniform distribution function.

\subsection{Series Systems}\label{sec:series}

This section considers the special case of problem $[KP_F]$ where the network is connected and each node contains at most one successor in addition to at most one predecessor, resulting in a series system. We assume that the CDF $F_Y(y)$  forms an S-curve such that $F_Y(y)$ is convex and nondecreasing for $y\le \beta $ and is concave and nondecreasing for $y\ge \beta$ for some $\beta\in [0,1]$.  Note that the extreme case of $\beta=0$ corresponds to a concave CDF, while the other extreme of $\beta=1$ corresponds to a convex CDF.  

For a series system, observe that any feasible solution that obeys Proposition \ref{p1} consists of an ordered set of $K$ groups such that $x_i=x_{\pi(i)}$ if both $i$ and its immediate predecessor are in the same group, while $x_i<x_{\pi(i)}$ if $i$ and its predecessor are in different groups. Consider two groups $G_1$ and $G_2$ such that $x_i=x_j$ for $(i,j) \in G_1$ and $x_i=x_j$ for $(i,j) \in G_2$, $x_i>x_j$ if $i \in G_1$ and $j \in G_2$. Let us assume that an optimal solution exists such that $\sum_{i=1}^n c_ix_i = B$ (instances that do not obey this are trivial because $F(\cdot)$ is nondecreasing, so that either an optimal solution exists satisfying this constraint at equality, or $\sum_{i=1}^n c_i < B$, in which case $x_i=1$ for all $i\in N$ in an optimal solution).  The following lemmas, along with Proposition \ref{p1}, allow us to consider a manageable subset of feasible solutions and lead to a polynomial-time algorithm for solving $[KP_F]$ when the network corresponds to a series system.

\begin{lemma}
    An optimal solution for a series system exists that contains at most one group $G_\cup$ with $x_i=x_j$ for $(i,j) \in G_\cup$ such that $x_i \in (0,\beta]$ for all $i \in G_\cup$.
    \label{l:gd_gt1}
\end{lemma}

\begin{lemma}
    An optimal solution for a series system exists that contains at most one group $G_\cap$ with $x_i=x_j$ for $(i,j) \in G_\cap$ such that $x_i \in (\beta, 1)$ for all $i\in G_{\cap}$.
    \label{l:gd_gt2}
\end{lemma}

We will utilize Lemmas \ref{l:gd_gt1} and \ref{l:gd_gt2} along with Proposition \ref{p1} to enumerate a polynomial number of solutions that satisfy these two lemmas, the proposition, and the problem's necessary KKT conditions.  

Observe that any feasible solution must consist of four sets $I_1=\{i: x_i=1\},\, I_0=\{i:x_i=0\},\, I_\cup = \{i:0<x_i\leq \beta\}$ and $I_\cap=\{i:\beta < x_i < 1\}$. We use the notation $I_1 \succ I_\cap \succ I_\cup \succ I_0$ to mean that nodes in $I_1$ must have $x_i$ values higher than any nodes in $I_\cap$, which must have $x_i$ values higher than any nodes in $I_\cup$, which must have $x_i>0$. Letting $I(i)$ denote the set that contains node $i$, Proposition \ref{p1} implies we can confine our search of candidate optimal solutions to those consisting of an allocation of indices to these four sets satisfying the following contiguity properties: 
\begin{enumerate}[(i)]
\item for all node pairs $i$ and $j$, where $i \in \mathcal{P}(j)$
, if $I(i)=I(j)$ and a directed path exists from $i$ to $j$, then all nodes on that directed path must also be in $I(i)$, \label{item1}
\item if a path exists from $i$ to $j$, then $I(i) \succeq I(j)$. \label{item2}
\end{enumerate}

In addition, for any allocation of indices to these subsets, only the values of $x_i$ for $i \in I_\cap \cup I_\cup$ are undetermined. For any such allocation, let $W_\cup=\sum_{i \in I_\cup} w_i,\, C_\cup=\sum_{i \in I_\cup} c_i,\,W_\cap=\sum_{i \in I_\cap} w_i,$ and $ C_\cap=\sum_{i \in I_\cap} c_i$, $x_\cup=x_i : i \in I_\cup$, and $x_\cap=x_i : i \in I_\cap$. 
For a given allocation of indices to subsets, applying Lemmas \ref{l:gd_gt1} and \ref{l:gd_gt2}, we wish to determine if a solution exists to the following system in the two variables $x_\cup$ and $x_\cap$:
\begin{eqnarray}
    & \frac{W_\cup}{C_\cup} f(x_\cup)  =\frac{W_\cap}{C_\cap}f(x_\cap), \nonumber \label{kkt1a} \\
    & C_\cup x_\cup + C_\cap x_\cap  = B-C_1,\\
    & 0 \leq x_\cup \leq \beta, \nonumber \\
    & \beta \leq x_\cap \leq 1, \nonumber \label{kkt4a}
\end{eqnarray}
where $C_1=\sum_{i \in I_1}c_i$.  This system can be solved using the following set of equations in a single variable by setting $x_\cup = \hat{B} - \hat{C}x_\cap$, where $\hat{B}=\frac{B-C_1}{C_\cup}$ and $\hat{C}=\frac{C_\cap}{C_\cup}$, and determining all solutions to the following system.
\begin{eqnarray}
    & f_Y(x_\cap)=\frac{W_\cup C_\cap}{W_\cap C_\cup} f_Y(\hat{B} - \hat{C}x_\cap), \label{kkt1b}  \\
    & \beta \leq x_\cap \leq 1. \label{kkt4b}
\end{eqnarray}
For a given allocation of items to the subsets $I_0,I_\cup,I_\cap,$ and $I_1$ satisfying contiguity properties (\ref{item1}) and (\ref{item2}), if a point can be found that satisfies \eqref{kkt1b}-\eqref{kkt4b}, then it is a KKT point. 
It is straightforward to show that for $n$ nodes in series, the number of candidate allocations of nodes to four subsets is $\mathcal{O}(n^3)$.  
If $\phi$ denotes the worst-case time required to enumerate all solutions to \eqref{kkt1b}-\eqref{kkt4b}, then the serial system case can be solved in $\mathcal{O}(\phi n^3)$.  The time required to enumerate solutions to \eqref{kkt1b}-\eqref{kkt4b} depends on the functional form of $f_Y(y)$. If, for example, $f_Y(y)$ corresponds to a triangular distribution, then $f_Y(y)$ is piecewise linear on two segments, and the unique solution to \eqref{kkt1b}-\eqref{kkt4b} (if one exists) can be determined in constant time. If, as in many commonly applied distributions, $f_Y(y)$ is unimodal with no flat intervals, then \eqref{kkt1b}-\eqref{kkt4b} contains at most one solution, which can be found (if one exists) using bisection search on the interval $[\beta,1]$.  The time required to obtain a solution within a tolerance of $\epsilon$ is $\mathcal{O}\left(\log_2 \left( \frac{1}{\epsilon}\right)\right)$.  Thus, for distributions including general versions of the Cauchy, gamma, Erlang, log-normal, logistic, normal, and Weibull distributions, the problem on a serial system can be solved in $\mathcal{O}\left( \log_2 \left( \frac{1}{\epsilon}\right)n^3 \right)$ time.  

Note that if the CDF is either convex or concave over the [0,1] interval, then only one of the sets $I_\cup$ or $I_\cap$ can exist. Moreover, in such cases, there are only three necessary subsets. If $F$ is concave (convex), then we have the sets $I_1$, $I_\cap$ ($I_\cup$), and $I_0$. For the concave case, \eqref{kkt1b}-\eqref{kkt4b} simply becomes 
\begin{center}
    $x_\cap=\frac{B-C_1}{C_\cap}$,\\
    $0 \leq x_\cap \leq 1$,
\end{center}
while if $F$ is convex, we obtain
\begin{center}
    $x_\cup=\frac{B-C_1}{C_\cup}$,\\
    $0 \leq x_\cup \leq 1$.
\end{center}
In both cases, this requires only determining whether the implied $x$ falls between 0 and 1 and the resulting worst-case complexity becomes $\mathcal{O}(n^2)$.  

We provide examples of the solutions for the system characterized by (\ref{kkt1b})-(\ref{kkt4b}) for some of the probability distributions (triangular, normal, uniform) that are commonly used to model the severity of disruptions \citep{azad2014new,dehghani2014impact} in Appendix~\ref{pdexamples}.

\subsection{Tree Networks with a Uniform Distribution}\label{sec:uniform}

Problem $[KP_F]$ when $F_Y(y)$ follows a uniform distribution on $[0,1]$ (so that $F_Y(y)=y$) becomes a linear program formulated as follows.
\begin{eqnarray*}
\mbox{Maximize} & \sum_{i=1}^n w_ix_i \label{m:gobju} \\
    \text{Subject to: } & \eqref{eq:capacity}-\eqref{gt_mod}. 
\end{eqnarray*}
We next provide a solution approach that enables faster solution of the problem in the worst case than is available using linear programming solvers.  Using the previously defined KKT multipliers as the dual variables corresponding to constraints \eqref{eq:capacity}-\eqref{eq:xi}, the corresponding dual problem for this LP can be written as
\begin{eqnarray}
\mbox{[D]} \quad \mbox{Minimize} & B \alpha + \lambda_1 & \\
\mbox{Subject to:} &   \lambda_i \ge \sum_{j\in \sigma(i)} \lambda_{j} + w_i - c_i \alpha & i \in N \backslash \mathcal{L}, \label{eq:treelambdai} \\
& \lambda_i \ge w_i - c_i \alpha & i \in\mathcal{L},\label{eq:treelambdan} \\
& \lambda_i \ge 0, & i = 1, \ldots, n, \\
& \alpha \ge 0. &
\end{eqnarray}
Observe that because we seek to minimize $\lambda_1$ for any $\alpha$, and $\lambda_1\ge \sum_{j\in \sigma(1)} \lambda_{j} + w_1 - c_1\alpha$, we can, without loss of optimality, set 
$$\lambda_1=\max\left\{\sum_{j\in \sigma(1)} \lambda_{j} + w_1 - c_1\alpha, 0\right\}=\left(\sum_{j\in \sigma(1)} \lambda_{j} + w_1 - c_1\alpha\right)^+.$$
Similarly, for any $j \in \sigma(1)$ we can set
$$\lambda_j=\left( \sum_{k\in \sigma(j)}\lambda_k+w_j-c_j\alpha\right)^+.$$
We can continue this substitution for successive $\lambda_k$ values, noting that for $i\in \mathcal{L}$ we obtain 
\begin{eqnarray*}
\lambda_i = \max \{w_i - c_i\alpha;0 \} = (w_i - c_i\alpha)^+.
\end{eqnarray*}
For any node $i$ such that all of its immediate successors are in $\mathcal{L}$, this will result in the expression
\begin{eqnarray*}
\lambda_{i} =  \left(\sum_{j\in \sigma(i)}(w_j - c_j\alpha)^+ + w_{i}-c_{i}\alpha\right)^+.
\end{eqnarray*}
In general, with a slight abuse of notation, we can write $\lambda_i$ as a recursive function of $\alpha$ corresponding to a summation of terms of the form 
\begin{eqnarray}
  \lambda_i(\alpha) = \left(\sum_{j\in \sigma(j)}\lambda_j(\alpha) +  w_i-c_i\alpha\right)^+, \label{eq:nu}
\end{eqnarray}
where each $\lambda_j(\alpha)$ takes the same form as in \eqref{eq:nu}, and $\sigma(j)=\emptyset$ for $j\in \mathcal{L}$.  We are interested in particular in characterizing the function $\lambda_1(\alpha)$, as we can rewrite the dual problem as  
\begin{eqnarray}
\mbox{[D']} \quad \mbox{Minimize} & z(\alpha)=B \alpha + \lambda_1(\alpha) & \\
\text{Subject to:} & \alpha \ge 0. &
\end{eqnarray}
\begin{lemma}\label{Lem:Piecewise}
$\lambda_1(\alpha)$ is a piecewise-linear, nonincreasing, and convex function of $\alpha$.
\end{lemma}
\begin{corollary}
The objective function of $[D']$, $z(\alpha)$, is a piecewise-linear and convex function of $\alpha$.
\end{corollary}
Because $z(0)=\sum_{i\in N} w_i$ and the slope of the piecewise-linear function $z(\alpha)$ is negative at $\alpha=0$ (with value $-\sum_{i \in N} c_i$), the optimal value of $\alpha$ for $[D']$ must be positive.  In addition, because $\lim\limits _{\alpha \rightarrow \infty}z(\alpha)=\infty$, the optimal value of $\alpha$ must be finite.  Thus, an optimal value of $\alpha$ must lie at a breakpoint of the function $z(\alpha)$.  It is possible to show that all breakpoints of $z(\alpha)$ can be expressed in the form $\frac{\sum_{i \in \mathcal{T}} w_i}{\sum_{i\in \mathcal{T}}c_i}$, for some subset $\mathcal{T}\subseteq N$ of the nodes in the tree network.  Moreover, if $\rho_{\min}=\min\limits_{i\in N}\left(\frac{w_i}{c_i}\right)$ and $\rho_{\max}=\max\limits_{i\in N}\left(\frac{w_i}{c_i}\right)$, then it is straightforward to show that
\begin{eqnarray}
    \rho_{\min} \le \frac{\sum_{i \in \mathcal{T}} w_i}{\sum_{i\in \mathcal{T}}c_i} \le \rho_{\max},
\end{eqnarray}
for any $\mathcal{T}\subseteq N$.  Because of this, an optimal value of $\alpha$ must lie on the interval $[\rho_{\min},\rho_{\max}]$. For any given $\alpha$, we can compute $z(\alpha)$ in $\mathcal{O}(n)$ time.  Because of this, letting $\ell=\rho_{\max}-\rho_{\min}$, we can perform bisection search over the interval $[\rho_{\min},\rho_{\max}]$ to determine an optimal solution (within a tolerance of $\epsilon$) in $\mathcal{O}\left( n \log\left(\frac{\ell}{\epsilon} \right)\right)$.

\section{Heuristic Solution Method for General Trees}\label{sec:heuristic}

We next propose a heuristic solution method for problem $[KP_F]$ for general tree systems under a disturbance severity distribution that takes the form of an S-shaped curve.  Recall that in the series system case, an optimal solution exists such that at most one contiguous group of nodes falls into each of the sets $I_1$, $I_\cap$, $I_\cup$, and $I_0$, i.e., the fortification level is equal to one, in the concave portion of the S-curve, in the convex portion of the S-curve, and equal to zero.  In the series system case, an optimal solution exists such that all nodes in $I_\cap$ have equal values, as do all nodes in $I_\cup$, which enabled deriving a polynomial-time solution algorithm.  In tree systems, although each node falls into one of these four sets, this property does not necessarily hold, i.e., nodes within each of the sets $I_\cap$ and $I_\cup$ do not necessarily have equal $x_i$ values.  

To find solutions for this more general class of problems, we propose a fast multi-start Network Search Algorithm that starts with an initial solution and makes local improvements at each iteration using procedures based on a priority metric that implicitly determines the `importance' of a node. Starting with multiple initial solutions typically enables quickly finding multiple locally optimal solutions. The starting solutions utilize the key properties of optimal solutions for series systems (in particular, Property \ref{p1} and Lemmas \ref{l:gd_gt1} and \ref{l:gd_gt2}). The initial solution (\ref{is:equal}) described below assumes all nodes have equal fortification levels, whereas (\ref{is:node1}) prioritizes the source node while maintaining equal levels for the remaining nodes. The initial solutions (\ref{is:concave}) and (\ref{is:convex}) below refer to node groups created based on the number of predecessors a node has ($|\mathcal{P}(i)|$, for each $i \in N$), i.e., nodes with the same number of predecessors are grouped together. These groups are then sorted in increasing order, with the first group consisting of only the source node (which has no predecessor). For the general tree example shown in Figure~\ref{fig:treegraph}, the aforementioned ordered groups of nodes would be formed as follows: $\{1\},\{2,3,5\},\{4,6,7,8\},\{9\},\{10,11,12\}, \{13\}$, and $\{14\}$. The initial solution (\ref{is:inbound}) uses a metric based on the weights and costs of the nodes on the path from the source to a particular node. An additional starting solution (\ref{envelope}) is obtained by using the solution that results when $F_Y(y)$ is replaced by its upper concave envelope, resulting in a convex program, the solution of which is feasible for problem $[KP_F]$.  


The set of initial feasible solutions consists of solutions constructed as follows: 
\begin{enumerate}[(i)]
\item \textbf{Equal Fortification:} All nodes are equally fortified at a value of $\frac{B}{\sum_{i=1}^n c_i}$; \label{is:equal}
\item \textbf{Fully Fortified Source:} Node 1 (the source) is fortified at a value of $x_1 =\min\left\{\frac{B}{c_1}, 1\right\}$ and all remaining nodes are equally fortified at a value of $\frac{B-c_1x_1}{\sum_{i=2}^n c_i}$; \label{is:node1}
\item \textbf{Depth-Based Concave:} Starting with the first group on the sorted list, nodes in the group are added to a set ($\mathcal{S}$) as long as equally fortified nodes in this set remain in the concave portion of the CDF. That is, node groups are added to $\mathcal{S}$ in sort order as long as $\beta \leq \frac{B}{\sum_{i \in \mathcal{S}} c_i} \leq 1$. We set the fortification level for each node $i\in \mathcal{S}$ equal to $\frac{B}{\sum_{i \in \mathcal{S}} c_i}$; \label{is:concave}

\item \textbf{Depth-Based Convex:} Starting with the first group on the sorted list, nodes in the group are added to a set ($\mathcal{S}$) until equally fortified nodes in this set remain in the convex portion of the CDF. That is, node groups are added to $\mathcal{S}$ in sort order as long as $0 \leq \frac{B}{\sum_{i \in \mathcal{S}} c_i} \leq \beta$. We set the fortification level for each node $i\in \mathcal{S}$ equal to $\frac{B}{\sum_{i \in \mathcal{S}} c_i}$; \label{is:convex}

\item \textbf{Cumulative Knapsack Ratio:} For each node, we compute the sum of the weights of all nodes on the path from the source and divide by the sum of the corresponding costs, i.e., $\phi_i=\frac{\sum_{j \in \mathcal{P}(i)} w_j}{\sum_{j \in \mathcal{P}(i)} c_j}$.  Defining $\mathcal{S}$ as the set of all nodes $i$ such that $\phi_i > \phi_{\pi(i)}$, plus all nodes $j\in \mathcal{P}(i)$, we set the fortification level for each node $i\in \mathcal{S}$ equal to $\frac{B}{\sum_{i \in \mathcal{S}} c_i}$;\label{is:inbound}
\item \textbf{Upper Concave Envelope Relaxation:} Nodes are fortified based on the optimal solution of the relaxation of $[KP_F]$ in which the CDF $F_Y(y)$ is replaced by its upper concave envelope. \label{envelope}
\end{enumerate}
 
Next, we describe the relaxation problem used to obtain the solution type (\ref{envelope}) above. Recall that we assumed that the disturbance severity distribution $F_Y(y)$ is differentiable and takes the form of an S-shaped curve, such that $F_Y(y)$ is convex and nondecreasing for $0 \le y\le \beta$ and is concave and nondecreasing for $\beta \le y \le 1$.  Assuming we can characterize the pdf $f_Y(y)$ (i.e., the derivative of $F_Y(y)$), we create the convex envelope for the S-curve by determining the smallest value of $y$ on the concave portion of the curve such that $\frac{F_Y(y)}{y}=F_Y'(y)=f_Y(y)$.  This corresponds to the point at which a line extending from the origin to the concave portion of the S-curve forms the concave envelope for the S-curve.  Thus, if $\tilde{y}$ solves $\min\left\{y\in [\beta,1]:F_Y(y)=yf_Y(y)\right\}$, then the concave envelope of $F_Y(y)$, which we denote as $\tilde{F}(y)$ is described as
\begin{eqnarray}
    \tilde{F}(y) = \left\{ \begin{matrix} \frac{f_Y(\tilde{y})}{\tilde{y}}y, & 0\le y\le \tilde{y}, \\ F_Y(y), & \tilde{y}\le y\le 1.\end{matrix} \right.
\end{eqnarray}

\begin{figure}[!ht]
    \centering
    \includegraphics[scale=0.5]{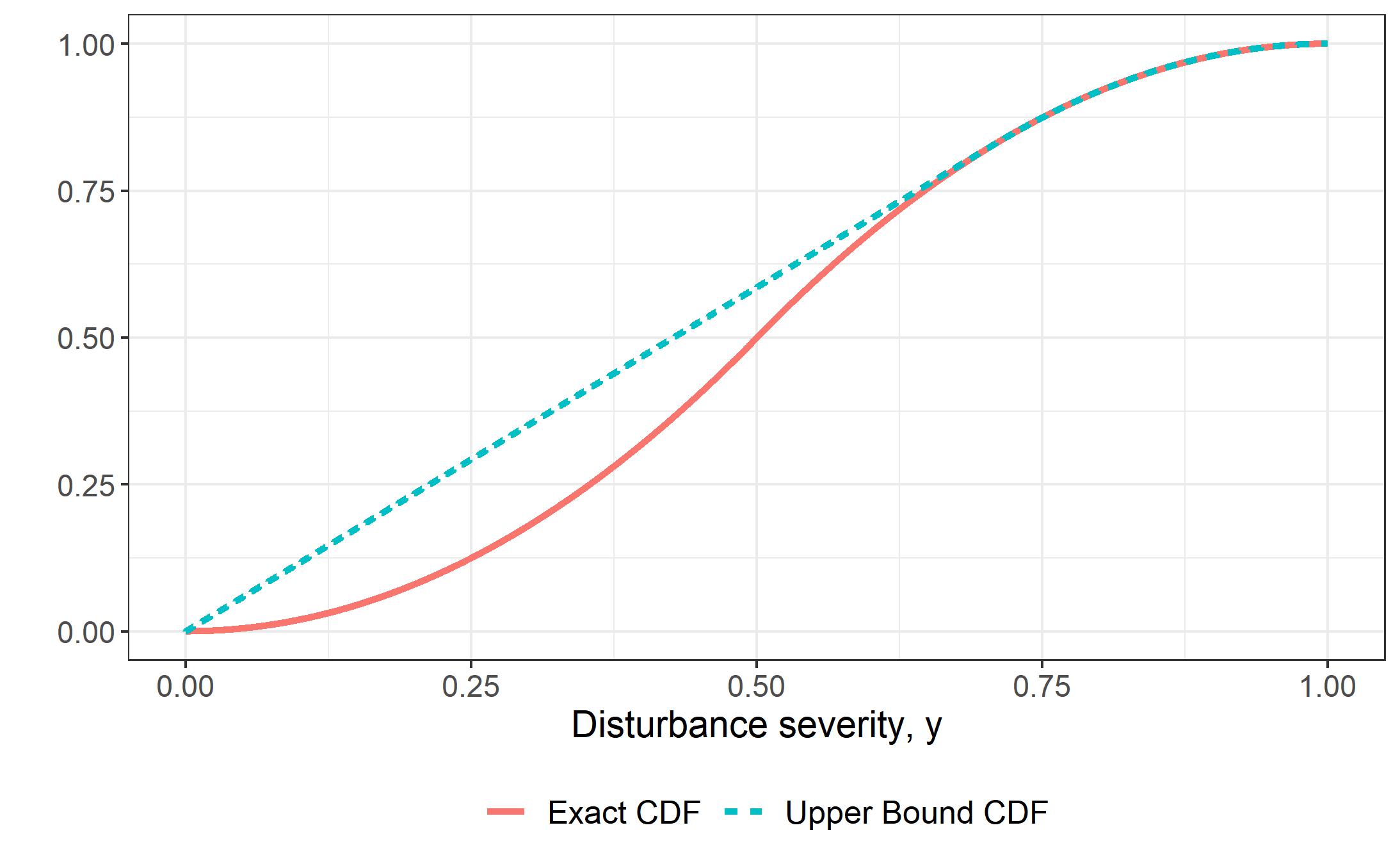}
    \caption{Upper concave envelope of the triangular distribution with $\beta=\frac{1}{2}$.} \label{fig:upper}
\end{figure}

If we replace each $F_Y(x_i)$ in $[KP_F]$ with $\tilde{F}(x_i)$, then the resulting problem is a convex program that can be solved quickly to obtain a candidate heuristic solution, and an initial solution for the Network Search Algorithm. Figure~\ref{fig:upper} illustrates the CDF and upper concave envelope of the triangular distribution when $\beta =\frac{1}{2}$. 

The existence of initial solutions (\ref{is:node1}), (\ref{is:concave}) and (\ref{is:convex}) is determined by the budget; if no such set of nodes can be formed, these solutions are omitted and not used as inputs to the algorithm. Note that with the exception of solution type (\ref{envelope}), all nodes with positive fortification levels that are less than one in the above initial solutions have equal fortification levels.  As a result, at least one of the sets $I_\cup$ and $I_\cap$ is empty for these solutions. 

Once we have generated all initial solutions, our Network Search Algorithm starts with an initial solution and makes local adjustments to node fortification values. Nodes are selected for adjustments in the Network Search Algorithm based on a priority measure (denoted by $\theta_i$) intended to implicitly reflect the relative `importance' of a node. This measure is associated with the marginal contribution of a given node $i$ to the objective function when $x_i=1$.  Although the initial form of this measure corresponds to the usual knapsack ratio, i.e., $\theta_i=\frac{w_i}{c_i}$, we consider alternative  measures of node importance, including $\theta_i=\frac{\sum_{j:i \in \mathcal{P}(j)} w_j}{c_i}$ and $\theta_i=\frac{\sum_{j\,:\, i \in \mathcal{P}(j)} w_j}{\sum_{j \in \mathcal{P}(i)} c_j}$ to account for the potential contribution of its successors and the cost of fortifying its predecessors, respectively. The different priority measures used in the Network Search Algorithm attempt to account for the way in which fortifying a node impacts other nodes in the network, while enabling fast identification of multiple candidate solutions.

Given an initial solution, the Network Search Algorithm applies four different solution adjustment procedures in an attempt to improve the existing solution. Each procedure is reapplied until a stopping criterion is met. That is, each procedure repeats until reaching a predetermined number of iterations without improvement before moving on to the next procedure. Each procedure selects a node for potential fortification-level modification from a given set of nodes based on the structure of the network. Examples of such sets are the set of leaf nodes among nodes with positive fortification levels  ($\mathcal{L}$) and the set of immediate successors of nodes in $\mathcal{L}$, denoted by $\Sigma(\mathcal{L})$ (For the general tree example shown in Figure~\ref{fig:treegraph}, these sets correspond to: $\mathcal{L}=\{4,6,7,10,11,14\}$, and $\Sigma(\mathcal{L})=\{3,5,9,13\}$). The \textit{Delete Node Procedure} chooses the lowest priority leaf node with a positive fortification level and sets its fortification level to zero. This frees up a portion of the limited budget and permits reallocating this budget to more critical nodes in the network. The \textit{Add Node Procedure} chooses the highest priority successor of a leaf node (whose fortification level is set to zero) and increases its fortification level. Upon completion of these two steps, the original network is typically reduced to a smaller set of nodes with the potential to increase the objective function value. The \textit{Max Fortification Procedure} then identifies the highest priority node among nodes with the highest degree in the reduced network. If this node's fortification level (and that of its predecessors) can be increased to one (by reallocating budget allocated to nodes at a value less than one), then this change is implemented. This increase in fortification levels requires reducing the budget allocated to those remaining nodes with fortification levels less than one in order to maintain feasibility. The $\text{lex}\!\max_{j \in N}\{\deg(j), \theta_j\}$ operator signifies a lexicographical preference for the maximization operation; that is, among the nodes with the highest degree, the highest priority node is chosen. 

Depending on whether the initial solution allocates node fortification levels in the convex region of the probability distribution, the concave region, or both (which applies to initial solution (\ref{envelope})), the \textit{Convex-Concave Set Assignment Procedure} uses three variations. In the first case, the node with the highest priority value among the highest degree nodes in the convex region is moved to the concave set ($I_\cap$) while updating the fortification levels of the remaining nodes accordingly. In the second case, the leaf node with the lowest priority in the concave region is added to the convex set ($I_\cup$), while in the last case, the highest priority node that was equal to zero is moved to the convex set ($I_\cup$), as are any of its predecessors that had a value of zero.  Any of its predecessors with a fortification level in the convex portion of the S-curve are then moved to the concave set ($I_\cap$). We provide a pseudocode for the Network Search Algorithm in Appendix~\ref{s:pseudo}.


After each procedure is completed, a new solution is obtained by solving the system of equations characterized by (\ref{kkt1a}). If this system has no solution, then we revert to the best solution obtained so far. After each procedure of the algorithm is applied, the difference between the objective function value of the new solution and that of the previous solution (denoted as $\Delta z$) is recorded, and the algorithm stops when the absolute value of $\Delta z$ is less than or equal to a pre-determined value ($\epsilon$). Note that the algorithm proceeds even if $\Delta z$ takes negative values, that is, the new solution is worse than the previous one for some maximum number of iterations, as the solution is not guaranteed to improve at every iteration. Allowing the exploration of seemingly worse solution directions often helps escape local optima and eventually find a better solution. Applying the Network Search Algorithm and using various definitions of the priority measure $\theta_i$ permits obtaining numerous candidate solutions, among which we select the one with the highest objective function value. After such a solution is selected, the process is then repeated for all initial solutions described above and the one that yields the greatest objective function value is returned as an output of the Network Search Algorithm.

\section{Computational Study}\label{s:computational}

In this section, we compare the results of the Network Search Algorithm (NSA) with the best available solution provided by a commercial global optimization solver. The mathematical model characterized by (\ref{gt_obj})-(\ref{gt_mod}) was coded using AMPL IDE (Version  4.0.1.202411072004), which used the BARON global optimization solver to obtain the best available solution for all randomly generated instances in the computational study within the time limit of 1 hour per problem instance. The NSA was coded using R Studio (Version 2022.07.2), which was also utilized to generate the random general tree instances. The initial solution (\ref{envelope}), Upper Concave Envelope, for the NSA was obtained by solving the relaxation of the problem $[KP_F]$ using the CONOPT nonlinear optimization solver in AMPL IDE. All elements of the computational study were conducted under the same hardware and software specifications (Windows 10, 12th Gen Intel(R) Core(TM) i7-12700 2.10 GHz, 16.0 GB RAM, 64-bit operating system, x64-based processor).

Our first set of computational experiments examines the performance of the NSA heuristic with respect to a variety of factors. Table~\ref{t:comp_par} describes the computational setting for randomly generated problem instances where the severity of the disruption follows a triangular (0,$\beta$,1) distribution for a general tree network. The triangular distribution is often used in the absence of information on the shape of the distribution or when there is little data available \citep{thomopoulos2017triangular}. The budget is varied as a function of the total cost to test the algorithm performance under low, medium, and high budget availability. Similarly, the maximum depth of the randomly generated network was varied (across the three values, 5, 10, and 20) to observe the impact of different structures on the solution, where a larger depth indicates a smaller number of branches exiting the source node.  We also considered how the value of $\beta$ affects performance by considering two uniform distributions for its random generation, $[0,1]$ and $[0,0.25]$.

\begin{table}[h]
    \caption{Computational Study Parameters}
    \centering
    \begin{tabular}{|rl|}
        \hline
        Parameter & Value\\
        \hline
         Number of nodes & $n$ = 50\\
        Node fortification costs & $c_i \sim$ Uniform$[1,10]$\\
        Node weights & $w_i \sim$ Uniform$[1,10]$\\
        Budget & $B \in$ $\left\{\frac{1}{5}\sum_{i=1}^n c_i,\frac{1}{2}\sum_{i=1}^n c_i,\frac{3}{4}\sum_{i=1}^n c_i\right\}$\\
        Network depth & $ \{5, 10, 20\}$\\
        Disruption severity  & $y\sim$ triangular(0,$\beta$,1)\\
        Distribution mode & $\beta \sim$ Uniform$[0,1]$ and $\beta \sim$Uniform$[0,0.25]$\\
        \hline
    \end{tabular}
    \label{t:comp_par}
\end{table}

Table~\ref{t:tri1} presents the average percentage difference between the objective function values of the algorithm and the best available solution provided by BARON for a variety of budget and network depth settings, where the severity of the disruption follows a triangular(0,$\beta$,1) distribution. Table~\ref{t:left} corresponds to a distribution mode $\beta$ randomly drawn from the [0,1] interval, while Table~\ref{t:right} corresponds to the interval [0,0.25]. The average percentage difference is defined as $\Delta \Bar{z}\% =  100\times\frac{z_\text{BARON}-z_\text{NSA}}{z_\text{BARON}}\%$, where $z_\text{BARON}$ and $z_\text{NSA}$ correspond to the objective function values of the best available solution provided by the global optimizer (BARON) and the NSA, respectively. Each percentage difference reported in Table~\ref{t:tri1} represents the average value among 30 instances of a randomly generated network with 50 nodes in the form of a general tree, using a computation time limit of 1 hour for each instance. Table~\ref{t:tri1} shows that as the budget increases, the gap between the NSA and solver solution decreases for both choices of the mode of distribution, $\beta$. Note that when the mode of distribution is chosen closer to the lower bound (as shown on Table~\ref{t:right}), which is more likely to be associated with real-life disasters where more severe events are less likely \citep{abedi2019review}, the overall NSA performance increases, with the average percentage difference decreasing from 0.4$\%$ to 0.1$\%$. A similar effect can be observed in the percentage of instances for which the heuristic algorithm is able to find the same solution as the global optimizer, which increases from 27$\%$ to 42$\%$ when the mode is closer to the lower bound. We observed a similar trend when we considered problem instances in which node fortification costs and node weights were drawn from a tighter interval [3,7] instead of [1, 10], i.e., when the variances of costs and weights are smaller, with every other parameter following the settings described in Table~\ref{t:comp_par}.  That is, the relative performance of the NSA heuristic improves with less variance in weights and costs, all else being equal.  We compare the objective function values of the algorithm and the solver (BARON) in Table~\ref{t:objvalues} of Appendix~\ref{s:objvalues}.

\begin{table}[h]
\centering
 \caption{$\Delta \Bar{z}\% $ with $n=50$ and triangular (0,$\beta$,1) distribution.}
\begin{minipage}{0.45\textwidth}
\addtocounter{table}{-1}
\renewcommand{\thetable}{\arabic{table}a}
\begin{tabular}{lrccc}
\cline{3-5}
\parbox[t]{1mm}{\multirow{3}{*}{\rotatebox[origin=c]{90}{\textbf{Budget}}}}      & High   &   \multicolumn{1}{|c}{0.1\% }   & 0.1\%          &  \multicolumn{1}{c|}{0.0\% }   \\
 & Medium & \multicolumn{1}{|c}{0.2\%}     & 0.1\%          & \multicolumn{1}{c|}{0.2\% }     \\
       & Low    & \multicolumn{1}{|c}{0.5\% }    & 0.9\%          & \multicolumn{1}{c|}{ 1.5\% }    \\
\cline{3-5}
                &        & Low       & Medium         & High      \\
                
                &        &  & \textbf{Depth} & 
\end{tabular} 

\caption{$\beta \sim$ Uniform$[0,1]$}\label{t:left}
\end{minipage}\hfill
\begin{minipage}{0.45\textwidth}
\addtocounter{table}{-1}
\renewcommand{\thetable}{\arabic{table}b}
\begin{tabular}{lrccc}
\cline{3-5}
\parbox[t]{1mm}{\multirow{3}{*}{\rotatebox[origin=c]{90}{\textbf{Budget}}}}      & High   &   \multicolumn{1}{|c}{0.0\% }   & 0.0\%          &  \multicolumn{1}{c|}{0.0\% }   \\
 & Medium & \multicolumn{1}{|c}{0.1\%}     & 0.0\%          & \multicolumn{1}{c|}{0.0\% }     \\
       & Low    & \multicolumn{1}{|c}{0.2\% }    & 0.3\%          & \multicolumn{1}{c|}{ 0.0\% }    \\
\cline{3-5}
                &        & Low       & Medium         & High      \\
                
                &        &  & \textbf{Depth} & 
\end{tabular}

\caption{$\beta \sim$Uniform$[0,0.25]$}\label{t:right}
\end{minipage}
    \label{t:tri1}
\end{table}

\begin{figure}[!tbp]
  \centering
  \subfloat[For varying maximum network depth]{\includegraphics[width=0.5\textwidth]{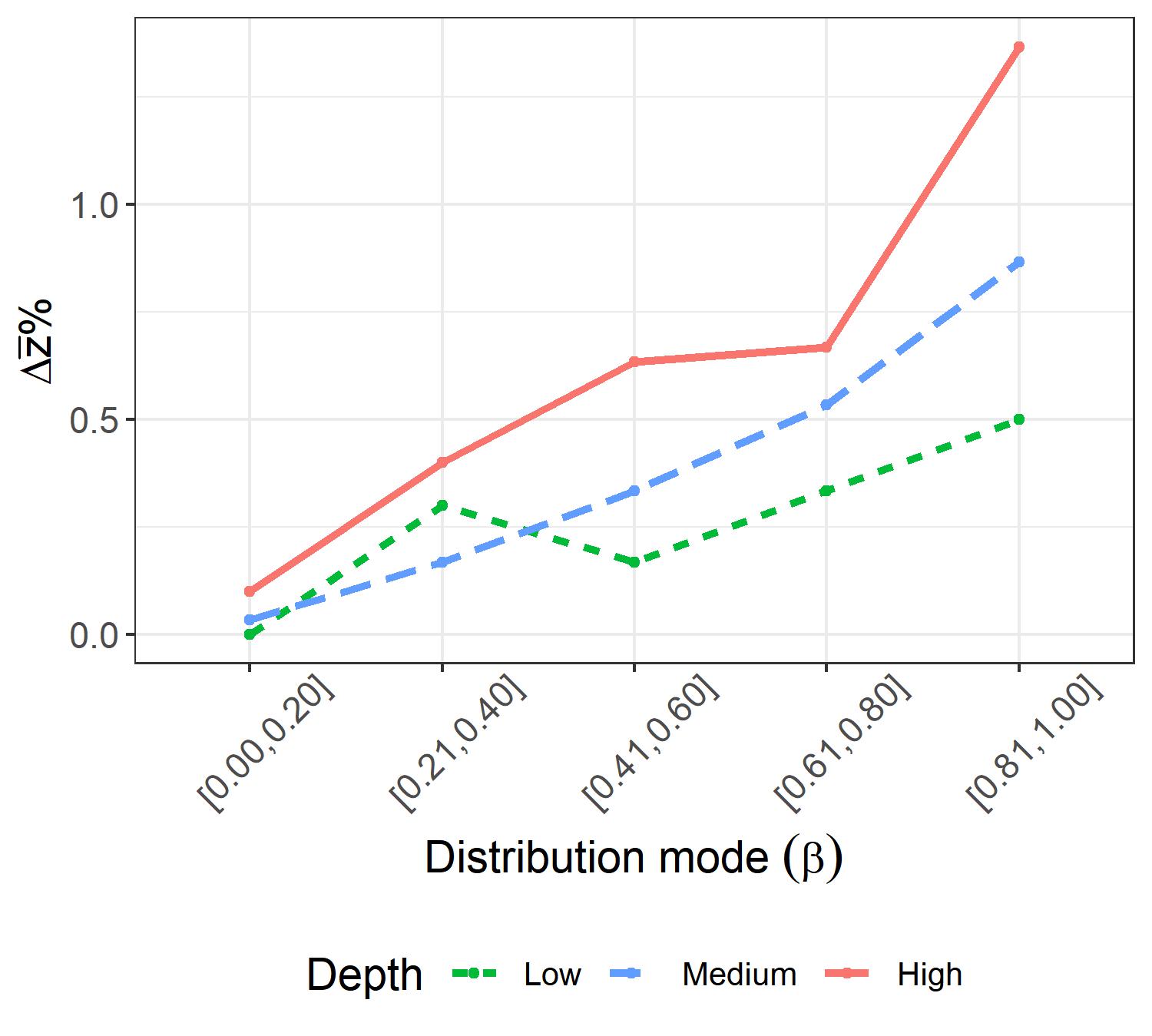}\label{fig:left}}
  \hfill
  \subfloat[For varying available budget]{\includegraphics[width=0.5\textwidth]{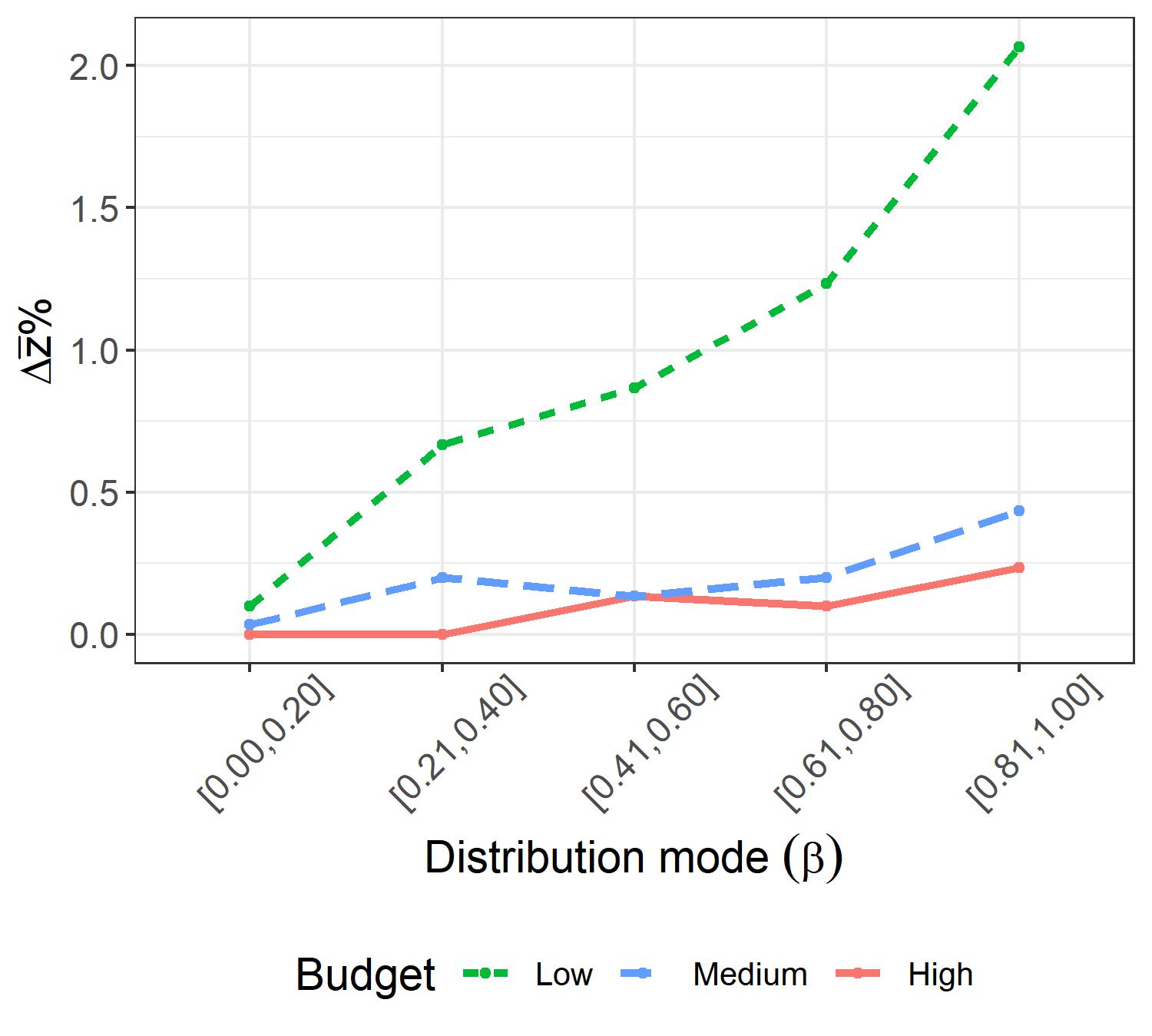}\label{fig:right}}
  \caption{$\Delta \Bar{z}\% $ with $n=50$ and triangular (0,$\beta$,1) distribution.} \label{fig:mode}
\end{figure}

Based on the problem instances in which $\beta \in [0,1]$, Figure~\ref{fig:mode} illustrates the change in the average percentage difference between the objective function values of the NSA and BARON for different values of the mode of the distribution, $\beta$, and for varying levels of maximum network depth and budget, respectively. Figure~\ref{fig:left} shows that in all but one of the intervals of $\beta$, the average percentage difference decreases as the network depth decreases, that is, as the number of branches in the network increases. Figure~\ref{fig:right} indicates that the NSA performs better when the available budget increases from low to medium. However, we do not observe such a steep change in the average percentage difference when the budget level rises from medium to high.

Figure~\ref{fig:left} shows that the average difference tends to decrease as $\beta$ decreases for all choices of maximum network depth. A similar observation can be made for all choices of budget in Figure~\ref{fig:right}.
This results because a smaller value of $\beta$ increases the concave portion of the CDF, and the initial solution (\ref{envelope}), Upper Concave Envelope, is obtained by solving a closer approximation of the CDF, i.e., the gap between the concave envelope and the actual CDF becomes smaller. That is, because the initial solution (\ref{envelope}) replaces the CDF with its upper concave envelope, this provides a near-optimal starting point for the algorithm at smaller values of $\beta$. Consequently, around 83$\%$ of the solutions provided by the NSA stem from the initial solution type (\ref{envelope}).

We further considered the impact of network size on algorithm performance by considering $n=100$ and $n=150$ node instances.  Table~\ref{t:times} presents the average computation times of the global optimizer (BARON) and our method, which consists of the sum of the computation times of the use of CONOPT to obtain initial solution (\ref{envelope}) and the NSA, for randomly generated tree networks with different numbers of nodes. The average computation times increase for both approaches as the number of nodes in the network increases. While the average computation time remains under a second for the NSA, the global optimizer hits the computational time limit as the number of nodes increases. For randomly generated general tree networks with 50 nodes, in approximately 29$\%$ of the cases, the global optimizer was unable to find an optimal solution within the given time limit. When the number of nodes in the network increases to 100, the solver reaches the time limit of an hour in all of the instances and is unable to find an optimal solution within the given time limit.

\begin{table}[h]
\caption{Average computation times for randomly generated general tree networks with a triangular(0,$\beta$,1) distribution and time limit of 3600 seconds.}
    \centering
    \begin{tabular}{|cc|c|c|c|}
    \hline
        \multicolumn{2}{|c|}{Approach}  & Average Time (sec) & Number of nodes & $\Delta \Bar{z} \%$\\
        \hline
        \multicolumn{2}{|c|}{BARON}  &1447.06  & \multirow{3}{*}{50} &\\
        \cline{1-3}
        \multirow{2}{*}{Our Method} & (NSA) & 0.13  &  & 0.4$\%$\\
         & (CONOPT) & 0.02 & &\\
        \hline
        \multicolumn{2}{|c|}{BARON} & 3600 & \multirow{3}{*}{100} & \\
        \cline{1-3}
        \multirow{2}{*}{Our Method}&(NSA) & 0.25 & & -5.5$\%$\\
        &(CONOPT) & 0.04 & &\\
        \hline
        \multicolumn{2}{|c|}{BARON} & 3600  & \multirow{3}{*}{150} & \\
        \cline{1-3}
        \multirow{2}{*}{Our Method}&(NSA) & 0.38 & &-3.3$\%$\\
        & (CONOPT) & 0.08 & &\\
        \hline
    \end{tabular}
    \label{t:times}
\end{table}

When comparing the best available solution provided by the global optimizer with the NSA solution for 100 nodes, the average percentage difference decreases to -5.5$\%$, which indicates that on average, the algorithm outperforms the global optimizer within the given time limit, with around 93$\%$ of the NSA solutions providing a greater objective function value than the solver. Similarly, when the number of nodes is increased to 150, BARON again reaches the time limit, and the average percentage difference between the objective functions of the algorithm and the global optimizer becomes -3.3$\%$, where the NSA solutions are strictly better than those obtained by BARON in all of the instances.

\section{Conclusion}\label{s:conclusion}

 We consider a distribution network for delivering a natural resource, service, or physical good to a set of nodes, each of which serves a set of customers, in which flow disruptions may occur at one or more nodes. We assume that all network nodes are subject to a common disturbance severity distribution, which is measurable and follows some well-defined probability distribution. For each node in the network, we wish to determine a fortification level that enables the node to withstand a disturbance up to some level of severity, given a limited fortification budget. The problem then seeks to minimize the expected number of customers without service.

 We formulate this problem as a continuous, nonlinear knapsack problem with precedence constraints and, although we show that the problem on a general tree is $\mathcal{NP}$-Hard, we characterize properties of optimal solutions that enable solving the problem for serial systems in polynomial time under mild conditions on the disturbance severity distribution. The difficulty of the resulting optimization problem depends both on the structure of the network and the corresponding disturbance severity distribution. For most commonly applied continuous distributions of natural phenomena, such as normal, beta, gamma, logistic, and triangular, the CDF is neither convex nor concave in general. For such cases, we develop an algorithm based on the problem's \emph{necessary} (but not sufficient) KKT conditions, which is able to eliminate the need for a global optimization solver, which is often unable to find an optimal solution in an acceptable amount of time as the problem size increases. The resulting algorithm is able to solve problems with a serial system structure in polynomial time under mild disturbance distribution assumptions. We use the properties of optimal solutions for series systems to develop a heuristic solution method, the Network Search Algorithm, for general tree distribution networks. We computationally test its performance against a commercial global optimization solver for a variety of randomly generated instances and show that the overall average percentage difference is less than $1\%$ while also maintaining the average computational time under a second. We also demonstrate the effectiveness of using the smallest concave envelope of the CDF to provide strong upper bounds and corresponding heuristic solutions. We note that enabling cycles in the distribution network is a valuable future direction for the continuation of our study, which can help extend the applicability of the findings to more systems. 

\bibliographystyle{chicago}
\spacingset{1}
\let\oldbibliography\thebibliography
 \renewcommand{\thebibliography}[1]{%
 	\oldbibliography{#1}%
 	\baselineskip11pt 
 	\setlength{\itemsep}{7pt}
 }
\bibliography{IISE-Trans}
\newpage
\appendix
\section{Appendix}

\label{App:Piecewise}
\subsection{Proofs of Lemmas and Propositions}\label{sec:Proofs}
\subsubsection{Proposition~\ref{p1}}

\begin{proof}
    Without loss of generality, let $x_j=x_i+\delta$ for some $i\in \mathcal{P}(j)$ where $\delta>0$ in a feasible solution for the network fortification problem. Then, $\min_{k \in \mathcal{P}(j)}\{x_k\} \neq x_j$, which indicates that the minimum value is attained at the fortification level of one of the predecessors of node $j$. Hence, letting $x_j'=x_j - \delta$ does not affect the probability $P\{y \geq \min_{k \in \mathcal{P}(j)}\{x_k\} \}$. Decreasing the fortification level of node $j$ in such a way generates extra budget of $\delta$ units in the problem without affecting the objective value with less cost. This shows that there exists another solution where $\delta$ can be used to increase the fortification level of predecessor nodes which would be at least as good as the current solution. Thus, in an optimal solution for the network fortification problem, $x_i \geq x_j$ for all $i\in \mathcal{P}(j)$.
\end{proof}

\subsubsection{Proposition~\ref{p:nphard}}
\begin{proof}
Consider a special case of $[KP_F]$ such that $c_i=w_i$, $i=1,\ldots,n$, with a source node 0 that supplies $n$ nodes indexed 1 through $n$ such that node 0 is the predecessor for each node $1,\ldots, n$.  Assume $F_Y(y)$ corresponds to any continuous distribution on $[0,1]$ such that $F_Y(y)\le y$ for $y\in [0, 1]$, with strict inequality holding for $y\in (0, 1)$.  Assume further that $c_0=0$ and $w_0>0$ so that $x_0=1$ in any optimal solution.  This special case of the problem is equivalent to the following optimization problem. 
\begin{eqnarray*}
        \mbox{[$KP_S$]} \qquad    \text{Maximize}  \quad  & \sum_{i=1}^n w_iF(x_i) \\
    \text{Subject to: } & \sum_{i=1}^n w_ix_i \leq B, \\
    & 0 \leq x_i \leq 1, & i=1,\ldots, n.
\end{eqnarray*}
The following S-curve distribution satisfies the stated conditions on $F_Y(y)$:
\begin{eqnarray*}
F_Y(y)=\left\{ \begin{matrix} k_1y^2, & 0\le y \le \beta, \\ y-k_2(1-y)^2, & \beta \le y\le 1, \end{matrix}   \right.
\end{eqnarray*}
where $k_1=\frac{1}{2\beta}+\frac{1}{2}$ and $k_2=\frac{\beta}{2(1-\beta)}$ for any $\beta\in (0,1)$. Figure \ref{fig:s-curve-dist} illustrates the shape of the pdf and CDF of the distribution when $\beta =\frac{2}{3}$.
\begin{figure}[!ht]
    \centering
    \includegraphics[scale=0.5]{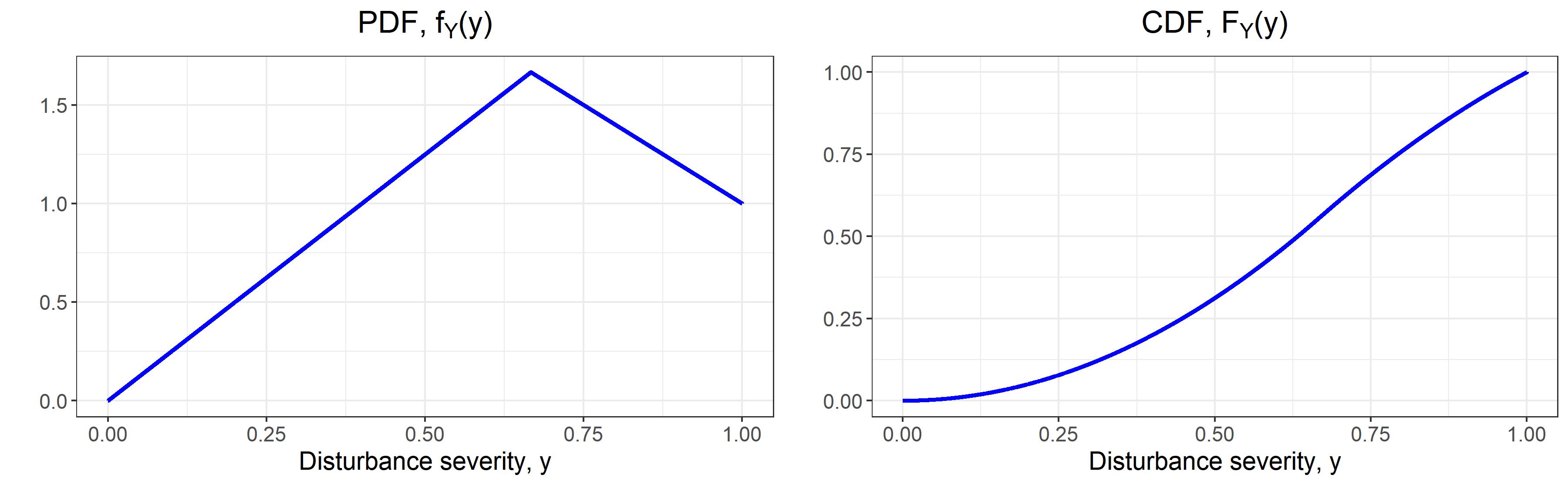}
    \caption{S-curve distribution with $\beta=\frac{2}{3}$.} \label{fig:s-curve-dist}
\end{figure}

\cite{kellerer2004subset} show that the following version of the subset sum problem is $\mathcal{NP}$-hard:
\begin{eqnarray*}
        \mbox{[$SSP$]} \qquad    \text{Maximize}  \quad  & \sum_{i=1}^n w_ix_i \\
    \text{Subject to: } & \sum_{i=1}^n w_ix_i \leq B, \\
    &  x_i \in \{0, 1\} & i=1,\ldots, n,
\end{eqnarray*}
where $w_i>0$ for $i=1,\ldots,n$.  We claim that [$KP_S$] has an optimal solution with objective function value $B$ if and only if [$SSP$] has an optimal solution with objective function value $B$, which occurs if and only if a subset of items $S\subset N$ exists such that $\sum_{i\in S} w_i=B$.  

To show this, observe that if a subset of items $S$ exists such that $\sum_{i\in S} w_i=B$, then by setting $x_i=1$ for $i\in S$, we obtain a solution for [$KPS$] with objective function value equal to $B$.  Such a solution must be optimal because $\sum_{i=1}^n w_i F(x_i) \le \sum_{i=1}^n w_i x_i \le B$.  Next, consider an optimal solution $x^*$ for [$KPS$] such that $\sum_{i=1}^n w_i F(x_i^*) = B$ and suppose $x_i^*\in (0,1)$ for at least one item $i$.  Because $F(x_i)<x_i$ for $x_i \in (0,1)$, it must be the case that $\sum_{i=1}^n w_i F(x_i^*)< \sum_{i=1}^n w_i x_i =B$, a contradiction, implying $x_i^*\in \{0, 1\}$ for all $i=1,\ldots,n$.  Because $\sum_{i=1}^n w_i F(x_i^*) = B$ and $\sum_{i=1}^n w_i F(x_i) \le \sum_{i=1}^n w_i x_i \le B$, it must be the case that $\sum_{i=1}^n w_i x_i^* = B$. Because $x_i^*\{0, 1\}$ for $i=1,\ldots,n$, this implies that a subset of items exists such that $\sum_{i\in S}w_i=B$.  Thus, any polynomial-time algorithm for solving $[KP_S]$ can solve $[SSP]$ in polynomial time, which can only hold if $\mathcal{P}=\mathcal{NP}$, implying that $[KP_S]$, and therefore $[KP_F]$, is $\mathcal{NP}$-hard.
\end{proof}

\subsubsection{Lemma~\ref{l:gd_gt1}}
\begin{proof}
    Suppose we have two groups of nodes, $G_1$ and $G_2$, with $x_i>x_j$ for $i \in G_1$ and $j \in G_2$ and both $0<x_j<x_i < \beta$. Suppose that $G_1$ and $G_2$ are chosen such that no other node $k$ satisfies $x_j<x_k<x_i$ (note that this is always possible if multiple groups exist with $x$ values on the convex portion of $F(x)$; moreover, the highest indexed node in $G_2$ must equal the highest indexed node in $G_1$ plus one). The capacity consumption by these two groups equals
    \begin{center}
        $x_iC(G_1)+x_jC(G_2),$
    \end{center}
    where $C(S)=\sum_{i \in S} c_i$. Suppose we increase $x_i$ by $\delta>0$ for $i \in G_1$. Then in order to maintain $\sum_{i=1}^n c_ix_i=B$ we decrease $x_j$ by $\delta \frac{C(G_1)}{C(G_2)}$. The new objective ($z_n$) minus the original one ($z_o$) is
    \begin{center}
        $z_n-z_o=(F(x_i+\delta)-F(x_i))W(G_1)+\left( F\left(x_j - \delta \frac{C(G_1)}{C(G_2)}\right)-F(x_j) \right)W(G_2)$,
    \end{center}
    where $W(S)=\sum_{i \in S} w_i$. If $F(x)$ is locally convex in the region containing $x_i$ and $x_j$, then for sufficiently small $\delta>0$
    \begin{center}
        $F(x_i+\delta) \geq F(x_i)+ \delta F'(x_i)$\\
        $F\left(x_j-\delta \frac{C(G_1)}{C(G_2)}\right) \geq F(x_j)-\delta\frac{C(G_1)}{C(G_2)}F'(x_j)$
    \end{center}
    This implies
    \begin{center}
        $z_n-z_o \geq \delta \left( F'(x_i)W(G_1)-\frac{C(G_1)}{C(G_2)} W(G_2)F'(x_j) \right)$
    \end{center}
    The new objective function is at least as high as the original one if
    \begin{center}
        $\frac{F'(x_i)}{F'(x_j)} \geq \frac{C(G_1)W(G_2)}{C(G_2)W(G_1)}$
    \end{center}
    If this holds then the original solution is either not optimal or we have an alternate optimal solution. We can continue this argument until one of the following occurs:
    \begin{itemize}
        \item The group containing $i$ is blocked from an increase by another group on the convex portion of $F(x)$, at which point it joins that group,
        \item The group containing $i$ hits $x_i=\beta$, at which point it becomes part of the concave portion of $F(x)$,
        \item The group containing $j$ is blocked from decrease by another group from below on the convex portion of $F(x)$, at which point it joins that group, or
        \item The group containing $j$ is blocked from decrease at $x_j=0$.
    \end{itemize}
    If, on the other hand,
    \begin{center}
        $\frac{F'(x_i)}{F'(x_j)} < \frac{C(G_1)W(G_2)}{C(G_2)W(G_1)}$
    \end{center}
    then we decrease $x_i$ by $\delta$ and increase $x_j$ by $\delta \frac{C(G_1)}{C(G_2)}$ and we have
    \begin{center}
        $F(x_i-\delta) \geq F(x_i)-\delta F'(x_i)$\\
        $F \left( x_j +\delta \frac{C(G_1)}{C(G_2)}  \right) \geq F(x_j) + \delta \frac{C(G_1)}{C(G_2)} F'(x_j)$
    \end{center}
    Thus,
    \begin{center}
        $z_n-z_o = W(G_1)(F(x_i-\delta)-F(x_i))+W(G_2) \left( F \left( x_j + \delta \frac{C(G_1)}{C(G_2)} \right) - F(x_j) \right)$\\
        $\geq \delta \left( W(G_2) \frac{C(G_1)}{C(G_2)} F'(x_j)-W(G_1)F'(x_i) \right) >0$
    \end{center}
    In this case, we can continue this until the group containing $x_i$ joins the group containing $x_j$. Repeated application of these arguments implies that at most one group can exist with equal $x_i$ values strictly between 0 and $\beta$, i.e., on the interior of the convex portion of $F(x)$.  Next, suppose we have an optimal solution such that $x_i=\beta$ for $i \in G_\beta$, as well as a group $G_\cup$ with $x_i=x_\cup$ for all $i \in G_\cup$ and $0<x_\cup <\beta<1$. Then we must have $\alpha=\frac{W(G_\beta)}{C(G_\beta)}F'(\beta)=\frac{W(G_\cup)}{C(G_\cup)}F'(\cup)$ at any associated KKT point. Suppose we increase $x_\cup$ by $\delta$ and decrease $x_\beta$ by $\delta\frac{C(G\cup)}{C(G_\beta)}$, leaving all other $x_i$ unchanged. Then the new solution minus the original one equals
    \begin{center}
        $z_n-z_o=W(G_\beta)\left( F\left(\beta-\delta\frac{C(G_\cup)}{C(G_\beta)}\right)-F(\beta) \right) + W(G_\cup)(F(x_\cup+\delta)-F(x_\cup))$
    \end{center}
    By convexity we have
    \begin{center}
        $F\left(\beta-\delta\frac{C(G_\cup)}{C(G_\beta)}\right) \geq F(\beta)-\delta\frac{C(G_\cup)}{C(G_\beta)} F'(\beta)$\\
        $F(x_\cup+\delta) \geq F(x_\cup)+\delta F'(x_\cup)$
    \end{center}
    so that
    \begin{center}
        $z_n-z_o \geq -W(G_\beta)\delta\frac{C(G_\cup)}{C(G_\beta)}F'(\beta)+W(G_\cup)\delta F'(x_\cup)$\\
        $=-W(G_\beta)\delta\frac{C(G_\cup)}{C(G_\beta)}\alpha\frac{C(G_\beta)}{W(G_\beta)}+W(G_\cup)\delta \alpha\frac{C(G_\cup)}{W(G_\cup)}$\\
        $=-C(G_\cup)\delta \alpha+ \alpha \delta C(G_\cup)=0.$
    \end{center}
    Thus, either $z_n=z_o$ and have an alternative optimal solution and can repeat this argument until $x_\beta=x_\cup$ or $z_n>z_o$ contradicting the optimality of the original solution.
\end{proof}

\subsubsection{Lemma~\ref{l:gd_gt2}}
\begin{proof}
    Suppose we have $G_1$ and $G_2$ defined as before, with $x_i>x_j$ for $i \in G_1$ and $j \in G_2$ and both $\beta < x_j<x_i < 1$. Suppose that $G_1$ and $G_2$ are chosen such that no other node $k$ satisfies $x_j<x_k<x_i$ (note that this is always possible if multiple groups exist with $x$ values on the concave portion of $F(x)$). If $F(x)$ is locally concave in the region containing $x_i$ and $x_j$, then for sufficiently small $\delta>0$
    \begin{center}
        $F(x_i) \leq F(x_i+\delta) -\delta F'(x_i + \delta)$\\
        $F(x_j) \leq F \left( x_j - \delta \frac{C(G_1)}{C(G_2)} \right) + \delta \frac{C(G_1)}{C(G_2)} F' \left( x_j -\delta \frac{C(G_1)}{C(G_2)} \right)$
    \end{center}
    Suppose we increase $x_i$ by $\delta$ and decrease $x_j$ by $\delta \frac{C(G_1)}{C(G_2)}$. Then
    \begin{center}
        $z_n-z_o=W(G_1)(F(x_i+\delta)-F(x_i))+W(G_2) \left( F \left( x_j - \delta \frac{C(G_1)}{C(G_2)} \right) -F(x_j) \right)$\\
        $\geq \delta \left( W(G_1) F'(x_i+\delta) - W(G_2)\frac{C(G_1)}{C(G_2)} F'\left( x_j - \delta \frac{C(G_1)}{C(G_2)} \right) \right)$
    \end{center}
    The new objective function is at least as high as the old one if
    \begin{center}
        $\frac{F'(x_i+\delta)}{F' \left( x_j - \delta \frac{C(G_1)}{C(G_2)}\right)}\geq \frac{C(G_1)W(G_2)}{W(G_1)C(G_2)}$
    \end{center}
    If this holds then the original solution is either not optimal or we have an alternate optimal solution.  We can continue this argument until one of the following occurs:
    \begin{itemize}
        \item The group containing $i$ is blocked from an increase by another group on the concave portion of $F(x)$, at which point it joins that group, 
        \item The group containing $i$ hits $x_i=1$,
        \item The group containing $x_j$ is blocked from decrease by another group from below on the concave portion of $F(x)$, at which point it joins that group, or
        \item The group containing $x_j$ decreases until $x_j=\beta$, at which point it is no longer on the interior of the concave portion of $F(x)$.
    \end{itemize}
    If
    \begin{center}
     $\frac{F'(x_i+\delta)}{F' \left( x_j - \delta \frac{C(G_1)}{C(G_2)}\right)} < \frac{C(G_1)W(G_2)}{W(G_1)C(G_2)}$
    \end{center}
    Then we decrease $x_i$ by $\delta$ and increase $x_j$ by $\delta \frac{C(G_1)}{C(G_2)}$. Then we get 
    \begin{center}
        $F(x_i) \leq F(x_i-\delta)+\delta F'(x_i-\delta)$\\
        $F(x_j) \leq F \left( x_j + \delta \frac{C(G_1)}{C(G_2)} \right) - \delta \frac{C(G_1)}{C(G_2)} F' \left( x_j + \delta \frac{C(G_1)}{C(G_2)} \right)$
    \end{center}
    Then
    \begin{center}
        $z_n-z_o = W(G_1)(F(x_i-\delta)-F(x_i))+ W(G_2) \left( F \left( x_j + \delta \frac{C(G_1)}{C(G_2)} \right) - F(x_j) \right)$\\
        $\geq W(G_2)\delta \frac{C(G_1)}{C(G_2)}F' \left( x_j + \delta \frac{C(G_1)}{C(G_2)} \right) - W(G_1) \delta F'(x_i -\delta)>0$
    \end{center}
    where the last inequality holds because $\frac{F'(x_i-\delta)}{F'\left(x_j +\delta \frac{C(G_1)}{C(G_2)} \right)} < \frac{F'(x_i + \delta)}{F' \left( x_j - \delta \frac{C(G_1)}{C(G_2)} \right)}$, which holds because $F(x)$ is nonnegative, concave, and nondecreasing and
    \begin{center}
     $\frac{F'(x_i+\delta)}{F' \left( x_j - \delta \frac{C(G_1)}{C(G_2)}\right)} < \frac{C(G_1)W(G_2)}{W(G_1)C(G_2)}$
    \end{center}
    by assumption. In this case, we can repeatedly apply this until the groups containing $i$ and $j$ join into one group. Repeated application of these arguments implies that at most one group can exist with equal $x_j$ values on the interval $(\beta, 1)$. 
\end{proof}

\subsubsection{Lemma~\ref{Lem:Piecewise}}
\begin{proof}
The function $\lambda_1(\alpha)$ is a nested collection of terms, with the innermost terms corresponding to leaf nodes (in the set $S$) with a value of $\max(w_i-c_i\alpha,0)$.  We start with these leaf nodes and work our way up the tree until we reach node 1, which corresponds to the outermost term.  For $i\in S$ and for $\alpha \le \frac{w_i}{c_i}$, $\lambda_i(\alpha) =\max(w_i-c_i\alpha,0)$ is a linear function with intercept $w_i$ and slope $-c_i$, and equals 0 for $\alpha\ge \frac{w_i}{c_i}$.  As the maximum between linear functions with non-positive slopes, $\max(w_i-c_i\alpha,0)$ is piecewise-linear, nonincreasing and convex in $\alpha$.  Each of these terms corresponding to a node in $S$ is embedded in a term associated with its predecessor. Next, we move up one level in the tree and consider any node $j$ whose only successors are in the set $S$.  For such a node, we have 
\begin{eqnarray}
\lambda_j(\alpha)=\left(\sum_{k\in \sigma(j)}\lambda_k(\alpha)+w_j-c_j\alpha\right)^+.   \label{eq:pcws}  
\end{eqnarray}
Because each $k\in S$, we established previously that the corresponding function $\lambda_k(\alpha)$ is a piecewise-linear, nonincreasing, and convex function.  Because the sum of piecewise-linear, nonincreasing, and convex functions is piecewise-linear, nonincreasing, and convex, the argument in parentheses in \eqref{eq:pcws} is a piecewise-linear, nonincreasing, and convex function.  The maximum between a nonincreasing, piecewise-linear, convex function and 0 is also a nonincreasing, piecewise-linear, convex function, implying that $\lambda_j(\alpha)$ is piecewise-linear, nonincreasing and convex. We can repeat this argument for each predecessor of a predecessor of a node in $S$ and obtain the same result.  Continuing up the tree until we reach node 1 implies the lemma.

\end{proof}

\subsection{Probability Distribution Examples} \label{pdexamples}
\noindent \textbf{Triangular Distribution}

A triangular distribution on [0,1] with the mode of distribution $\beta$, has 
\begin{center}
    $F'(x)=\begin{cases}
        \frac{2}{\beta}x, \;\;\;\;\;\;\;\;\;\;\;\;\;\;\;\; 0 \leq x \leq \beta\\
        \frac{2}{1-\beta}(1-x), \;\;\; \beta \leq x \leq 1
    \end{cases}$
\end{center}
In this case the system (\ref{kkt1a}) becomes 
\begin{center}
     $\frac{2W_\cup}{\beta C_\cup} x_\cup=\frac{W_\cap}{C_\cap}\frac{2}{1-\beta} (1-x_\cap)$\\
     $C_\cup x_\cup + C_\cap x_\cap  = B-C_1$\\
     $0 \leq x_\cup \leq \beta$\\
     $\beta \leq x_\cap \leq 1$
\end{center}
Thus it is a simple matter to determine whether this linear system has a solution in constant time. Solving for $x_\cup$ and $x_\cap$, we obtain the following
\begin{center}
    $x_\cup=\frac{C_\cup W_\cap \beta (B-C_1-C_\cap)}{C_\cup^2 W_\cap \beta - C_\cap^2 W_\cup (1-\beta)}$\\
    $x_\cap=\frac{B-C_1}{C_\cap}-\frac{C_\cup}{C_\cap}\frac{C_\cup W_\cap \beta (B-C_1-C_\cap)}{(C_\cup^2 W_\cap  + C_\cap^2 W_\cup )\beta -C_\cap^2 W_\cup}$
\end{center}

\noindent\textbf{Normal Distribution}

Suppose we have a normal distribution on the [0,1] interval such that $P(X<0)$ and $P(X>1)$ are negligible, with mean $\beta$ and variance $\sigma^2$. In this case, $F'(x)$ is the pdf of the normal distribution, i.e., 
\begin{center}
    $\frac{1}{\sqrt{2\pi}\sigma}e^{-\left( \frac{x-\beta}{\sigma}\right)^2}$
\end{center}
This gives the following system:
\begin{center}
     $\frac{W_\cup}{C_\cup} e^{-\left( \frac{x_\cup-\beta}{\sigma}\right)^2} =\frac{W_\cap}{C_\cap} e^{-\left( \frac{x_\cap-\beta}{\sigma}\right)^2}$  \\
     $C_\cup x_\cup + C_\cap x_\cap  = B-C_1$\\
     $0 \leq x_\cup \leq \beta$\\
     $\beta \leq x_\cap \leq 1$ 
\end{center}
Taking natural logs of both sides of the first equation gives
\begin{center}
    $(\beta-x_\cup)^2=(x_\cap-\beta)^2-\sigma^2 \left( \ln{ \left( \frac{W_\cap}{C_\cap} \right) } - \ln{ \left( \frac{W_\cup}{C_\cup} \right) } \right)$
\end{center}
Letting $\Lambda=\sigma^2 \left( \ln{ \left( \frac{W_\cap}{C_\cap} \right) } - \ln{ \left( \frac{W_\cup}{C_\cup} \right) } \right)$ and taking the square root of both sides gives
\begin{center}
    $x_\cup=\beta-\sqrt{(x_\cap-\beta)^2-\Lambda}$
\end{center}
Plugging this into the second equation gives
\begin{eqnarray}
    C_\cap x_\cap -C_\cup \sqrt{(x_\cap-\beta)^2-\Lambda}=B-C_1-C_\cup\beta
    \label{normaleq}
\end{eqnarray}
Note that if $\Lambda=0$ (if $\frac{W_\cap}{W_\cap}=\frac{W_\cup}{C_cup}$), then the left-hand side equals $C_\cap x_\cap-C_\cup(x_\cap-\beta)$, and then we get the solution
\begin{center}
    $x_\cap=\frac{B-C_1-2C_\cup \beta}{C_\cap - C_\cup}$
\end{center}
If $\lambda>0$, which occurs when $\frac{W_\cap}{W_\cap}>\frac{W_\cup}{C_\cup}$, then this system can only have a real solution for $x_\cap \geq \beta + \sqrt{\Lambda}$. Letting $g(x_\cap)=C_\cap x_\cap -C_\cup \sqrt{(x_\cap-\beta)^2-\Lambda}$, then \ref{normaleq} can be written as $g(x_\cap)=\alpha$, and we have $g'(x_\cap)=C_\cap  -C_\cup \frac{x_\cap - \beta}{\sqrt{(x_\cap-\beta)^2-\Lambda}}$. If $\Lambda=0$, then $g'(x_\cap)=C_\cap  -C_\cup$. If $\Lambda>0$ ($\frac{W_\cap}{W_\cap}>\frac{W_\cup}{C_\cup}$), then $g'(x_\cap)<C_\cap  -C_\cup$, and if $\Lambda<0$ ($\frac{W_\cap}{W_\cap}<\frac{W_\cup}{C_\cup}$), then $g'(x_\cap)>C_\cap  -C_\cup$. Observe that
\begin{center}
    $g''(x_\cap)=\frac{(x_\cap - \beta)^2}{((x_\cap-\beta)^2-\Lambda)^\frac{3}{2}}-\frac{1}{((x_\cap-\beta)^2-\Lambda)^\frac{1}{2}}$\\
    $=\frac{(x_\cap - \beta)^2}{((x_\cap-\beta)^2-\Lambda)^\frac{3}{2}}-\frac{(x_\cap-\beta)^2-\Lambda}{((x_\cap-\beta)^2-\Lambda)^\frac{3}{2}}$\\
    $\frac{\Lambda}{((x_\cap-\beta)^2-\Lambda)^\frac{3}{2}}$
\end{center}
Therefore, $g(x_\cap)$ is concave in $x_\cap$ if $\Lambda<0$ and is convex in $x_\cap$ if $\Lambda>0$. In either case, $g(x_\cap)=\alpha$ for any scalar $\alpha$ can have at most two solutions. If $\Lambda<0$ and $C_\cap>C_\cup$, then $g(x_\cap)$ is concave and increasing so at most one solution exists. Similarly, if $\Lambda>0$ and $C_\cap<C_\cup$, then $g(x_\cap)$ is convex and decreasing so at most one solution exists. In cases where $g(x_\cap)$ is not strictly increasing or decreasing, then we can determine $x_\cap$ such that $g'(x_\cap)=0$ and perform a gradient search (perhaps to the right and left of the point at which $g'(x_\cap)=0$) ) to determine whether a point exists such that $g(x_\cap)=\alpha$ for some $\beta \leq x_\cap \leq 1$. This should be able to be accomplished efficiently using bisection search. 


\noindent\textbf{Uniform Distribution}

A uniform [0,1] distribution can be considered as a special case of the differentiable distribution functions that are convex for $x \leq \beta$ and concave for $x \geq \beta$ by setting $\beta$ equal to either of the end points. This indicates that either $I_\cap$ or $I_\cup$ can exist, hence we only have three subsets $I_1,I_\cap$ ($I_\cup$) and $I_0$. The system (\ref{kkt1a}) becomes:

\begin{center}
    $  x_\cap   = \frac{B-C_1}{C_\cap}$\\
    $ 0 \leq x_\cap  \leq 1$\\
\end{center}
\subsection{The Pseudocode for the Network Search Algorithm} \label{s:pseudo}

\begin{algorithm}

\caption{Network Search Algorithm}
\label{alg-ns}
\begin{algorithmic}
\begin{singlespace}
    \State \textbf{Initialization:} Start with an initial solution with corresponding sets $I_1$, $I_\cap$, $I_\cup$, and $I_0$, and a priority metric $\theta$.
    \State Define $\mathcal{N}$ as the set of nodes with positive fortification levels, $\mathcal{L}$ as the initial set of leaf nodes in the subgraph induced by the nodes in $\mathcal{N}$, and $\mathcal{F}\subseteq\mathcal{N}$ as the set of nodes with positive fortification levels less than one.  Given a leaf node set $\mathcal{L}$, let $\Sigma(\mathcal{L})$ denote the set of immediate successors of nodes in $\mathcal{L}$.
    \If{$I_\cup \ne \emptyset$} 
        \State \textbf{Delete Node Procedure:}
        \State Choose $i=\arg\!\min_{j \in \mathcal{L}} \{\theta_j\}$, set $x_i=0$, and add $i$ to $I_0$
        \State $\mathcal{N}\gets \mathcal{N} \backslash\left\{i\right\}$, $\mathcal{F}\gets \mathcal{F} \backslash\left\{i\right\}$, $\mathcal{L}\gets\left(\mathcal{L}\backslash\{i\}\right)\cup\{\pi(i)\}$ and update $\Sigma({\mathcal{L}})$ 
        \State For each $i\in \mathcal{F}$, set $x_i=\min \left\{ \frac{B-\sum_{j\in I_1}c_j}{\sum_{j\in \mathcal{F}}c_j}, 1\right\}$  
        \State \textbf{if} $x_i=1$ for $i\in\mathcal{F}$ \textbf{then} $I_1\gets I_1\cup \mathcal{F}$ and $I_\cap=\emptyset$
        \State \textbf{else} $I_\cap=\mathcal{F}$ and $I_\cup=\emptyset$
        \State \textbf{Add Node Procedure:}
        \State Choose $i=\arg\!\max_{j \in \Sigma(\mathcal{L})} \{\theta_j\}$ 
        \State $\Sigma(\mathcal{L})\gets\left( \Sigma(\mathcal{L})\cup\{i\} \right)\backslash \pi(i)$, $\mathcal{F}\gets\mathcal{F}\cup\{i\}$, $\mathcal{N}\gets\mathcal{N}\cup\{i\}$, and $I_0\gets I_0\backslash\{i\}$
        \State For each $i\in \mathcal{F}$, set $x_i=\frac{B-\sum_{j\in I_1}c_j}{\sum_{j\in \mathcal{F}}c_j}$
        \State \textbf{if} $x_i<\beta$ for $i\in\mathcal{F}$ \textbf{then} $I_\cup=\mathcal{F}$ and $I_\cap=\emptyset$
        \State \textbf{else} $I_\cap=\mathcal{F}$ and $I_\cup=\emptyset$
        \State \textbf{Max Fortification Procedure:}
        \State Choose $i=\text{lex}\!\max_{j \in \mathcal{F}}\left\{\deg(j), \theta_j\right\}$, and let $\tilde{x}=\min\left\{\frac{B-\sum_{k\in I_1}c_k}{\sum_{k\in \mathcal{F}\cap\mathcal{P}(i)}c_k},1 \right\}$
        \State \textbf{if} $\tilde{x}=1$ \textbf{then} set $x_j=1$ for $j\in\mathcal{F}\cap\mathcal{P}(i)$, $I_1\gets I_1\cup\left(\mathcal{F}\cap\mathcal{P}(i)\right)$ and $\mathcal{F}\gets\mathcal{F}\backslash\left(\mathcal{F}\cap\mathcal{P}(i)\right)$
        \State \hspace{\algorithmicindent}  For each $i\in \mathcal{F}$, set $x_i=\frac{B-\sum_{j\in I_1}c_j}{\sum_{j\in \mathcal{F}}c_j}$ 
        \State \hspace{\algorithmicindent} \textbf{if} $x_i<\beta$ for $i\in\mathcal{F}$ \textbf{then} $I_\cup=\mathcal{F}$ and $I_\cap=\emptyset$
        \State \hspace{\algorithmicindent} \textbf{else} $I_\cap=\mathcal{F}$ and $I_\cup=\emptyset$
        \State \textbf{else} continue
        \State \textbf{Convex-Concave Set Assignment (I):}
        \State Choose $i=\text{lex}\!\max_{j \in I_\cup}\{\deg(j), \theta_j\}$, $I_\cap \gets I_\cap \cup \mathcal{P}(i)$ and $I_\cup \gets I_\cup \backslash \mathcal{P}(i) $
        \State Solve system (\ref{kkt1a}) based on the assignment of nodes to sets
    \ElsIf{$I_\cap \ne \emptyset$}
        \State Apply \textbf{Delete Node Procedure}
        \State Apply \textbf{Add Node Procedure}
        \State \textbf{Convex-Concave Set Assignment Procedure (II):}
        \State Choose $i=\arg\!\min_{j \in \mathcal{L}\cap I_\cap} \{\theta_j\}$, update $I_\cup \gets I_\cup \cup \{i\} $ and  $I_\cap \gets I_\cap \backslash \{i \} $
        \State Solve system (\ref{kkt1a}) based on the assignment of nodes to sets
        \State Apply \textbf{Max Fortification Procedure}
    \Else
        \State \textbf{Convex-Concave Set Assignment Procedure (III):}
        \State Choose $i=\arg\!\max_{j \in I_0}\{ \theta_j\}$, update $I_\cap \gets I_\cap \cup \{j\} $ for all $j \in \mathcal{P}(i) \cap I_\cup$  
        \State and $I_\cup \gets I_\cup \cup \{j\} $ for all $j \in \mathcal{P}(i) \cap I_0$
        \State Solve system (\ref{kkt1a}) based on the assignment of nodes to sets
    \EndIf
\end{singlespace}

\end{algorithmic}
\end{algorithm}

\newpage
\subsection{Computational Results for Selected Instances} \label{s:objvalues}

Table~\ref{t:objvalues} compares the objective function values and the percentage differences between these values for the Network Search Algorithm (NSA) and the global optimizer (BARON) for randomly generated general tree networks with 50 nodes. For each maximum network depth and budget configuration, 3 instances were chosen to represent the minimum, average, and maximum percentage differences observed respectively, out of 30 randomly generated instances. We observe that the maximum percentage difference decreases as the budget increases.
\begin{table}[h]
\caption{The objective function values of the NSA and the solver (BARON), and the average percentage differences between these values for randomly generated general tree networks with 50 nodes using triangular (0,$\beta$,1) distribution where $\beta$ is randomly generated between [0,1] for different levels of maximum network depth and budget.}
\centering
\begin{tabular}{|c|c|c|c|c|}
\hline
\multicolumn{2}{|c|}{Objective Function Values} & Percentage Differences & \multirow{3}{*}{Depth} & \multirow{3}{*}{Budget} \\
\cline{1-2}
\multirow{2}{*}{Algorithm} & Global Optimizer &  Between Objective  &   &  \\
& (BARON) & Function Values & & \\
\hline
108.40    & 108.33                   & -0.1\%                &     & \multirow{9}{*}{Low}    \\
112.55    & 113.15                   & 0.5\%                 & Low    &     \\
108.16    & 110.72                   & 2.3\%                 &     &     \\
\cline{1-4}
106.91    & 106.70                   & -0.2\%                &  &     \\
112.61    & 113.49                   & 0.8\%                 & Medium &     \\
71.32     & 77.01                    & 7.4\%                 &  &     \\
\cline{1-4}
108.82    & 108.02                   & -0.7\%                &    &     \\
100.37    & 101.86                   & 1.5\%                 & High   &    \\
67.24     & 73.87                    & 9.0\%                 &    &    \\
\hline
242.62    & 242.49                   & -0.1\%                &     & \multirow{9}{*}{Medium} \\
176.64    & 176.94                   & 0.2\%                 & Low    &  \\
211.80    & 214.47                   & 1.2\%                 &     &  \\
\cline{1-4}
225.68    & 224.87                   & -0.4\%                &  &  \\
180.26    & 180.36                   & 0.1\%                 & Medium &  \\
169.00    & 170.19                   & 0.7\%                 &  &  \\
\cline{1-4}
191.41    & 190.42                   & -0.5\%                &    &  \\
200.48    & 200.80                   & 0.2\%                 & High   &  \\
167.19    & 170.84                   & 2.1\%                 &    &  \\
\hline
268.60    & 268.39                   & -0.1\%                &     & \multirow{9}{*}{High}   \\
251.47    & 251.82                   & 0.1\%                 & Low    &    \\
279.01    & 280.57                   & 0.6\%                 &     &    \\
\cline{1-4}
250.47    & 250.43                   & 0.0\%                 &  &    \\
288.14    & 288.34                   & 0.1\%                 & Medium &    \\
221.20    & 222.95                   & 0.8\%                 &  &    \\
\cline{1-4}
284.75    & 283.36                   & -0.5\%                & \multirow{3}{*}{High}   &    \\
276.45    & 276.45                   & 0.0\%                 &    &    \\
235.52    & 237.89                   & 1.0\%                 &    &   \\
\hline
\end{tabular}
\label{t:objvalues}
\end{table}
\end{document}